\PassOptionsToPackage{unicode}{hyperref}
\PassOptionsToPackage{hyphens}{url}
\documentclass[
  11pt,
  a4paper,
]{article}
\usepackage{amsmath,amssymb}
\usepackage{iftex}
\ifPDFTeX
  \usepackage[T1]{fontenc}
  \usepackage[utf8]{inputenc}
  \usepackage{textcomp} % provide euro and other symbols
\else % if luatex or xetex
  \usepackage{unicode-math} % this also loads fontspec
  \defaultfontfeatures{Scale=MatchLowercase}
  \defaultfontfeatures[\rmfamily]{Ligatures=TeX,Scale=1}
\fi
\ifPDFTeX\else
\fi
\IfFileExists{upquote.sty}{\usepackage{upquote}}{}
\IfFileExists{microtype.sty}{% use microtype if available
  \usepackage[]{microtype}
  \UseMicrotypeSet[protrusion]{basicmath} % disable protrusion for tt fonts
}{}
\makeatletter
\@ifundefined{KOMAClassName}{% if non-KOMA class
  \IfFileExists{parskip.sty}{%
    \usepackage{parskip}
  }{% else
    \setlength{\parindent}{0pt}
    \setlength{\parskip}{6pt plus 2pt minus 1pt}}
}{% if KOMA class
  \KOMAoptions{parskip=half}}
\makeatother
\usepackage{xcolor}
\usepackage{longtable,booktabs,array}
\usepackage{calc} % for calculating minipage widths
\usepackage{etoolbox}
\makeatletter
\patchcmd\longtable{\par}{\if@noskipsec\mbox{}\fi\par}{}{}
\makeatother
\IfFileExists{footnotehyper.sty}{\usepackage{footnotehyper}}{\usepackage{footnote}}
\makesavenoteenv{longtable}
\usepackage{graphicx}
\makeatletter
\def\maxwidth{\ifdim\Gin@nat@width>\linewidth\linewidth\else\Gin@nat@width\fi}
\def\maxheight{\ifdim\Gin@nat@height>\textheight\textheight\else\Gin@nat@height\fi}
\makeatother
\setkeys{Gin}{width=\maxwidth,height=\maxheight,keepaspectratio}
\makeatletter
\def\fps@figure{htbp}
\makeatother
\providecommand{\tightlist}{%
  \setlength{\itemsep}{0pt}\setlength{\parskip}{0pt}}
\usepackage[a4paper,margin=2.5cm]{geometry}
\usepackage{booktabs,longtable,array,calc}
\usepackage{amsmath,amssymb}
\usepackage{graphicx}
\usepackage[hidelinks]{hyperref}
\usepackage{microtype}
\ifPDFTeX
  \usepackage{mathptmx}
  \DeclareUnicodeCharacter{2248}{\ensuremath{\approx}}
  \DeclareUnicodeCharacter{25A1}{\ensuremath{\square}}
\else
  \usepackage{fontspec}
  \usepackage{unicode-math}
\fi
\usepackage{setspace}
\usepackage{titling}\pretitle{\begin{center}\Large\bfseries}\posttitle{\end{center}}

\providecommand{\tightlist}{\setlength{\itemsep}{0pt}\setlength{\parskip}{0pt}}
\usepackage{enumitem}\setlist{nosep}
\usepackage{xurl}
\usepackage{etoolbox}\AtBeginEnvironment{longtable}{\small}
\ifLuaTeX
  \usepackage{selnolig}  % disable illegal ligatures
\fi
\IfFileExists{bookmark.sty}{\usepackage{bookmark}}{\usepackage{hyperref}}
\IfFileExists{xurl.sty}{\usepackage{xurl}}{} % add URL line breaks if available
\hypersetup{
  pdftitle={Time-reliability optimization of multi-pass underfill dispensing sequences on printed circuit boards: a schedule-dependent reliability model solved by a multiset-permutation binary-addition-tree algorithm},
  pdfauthor={Wei-Chang Yeh},
  hidelinks,
  pdfcreator={LaTeX via pandoc}}

\title{Time-reliability optimization of multi-pass underfill dispensing sequences on printed circuit boards: a schedule-dependent reliability model solved by a multiset-permutation binary-addition-tree algorithm}
\author{Wei-Chang Yeh}
\date{Integrated \& Collaborative AI Robotics Laboratory (ICARL), Department of Industrial Engineering and Engineering Management, National Tsing Hua University, Hsinchu, Taiwan. Corresponding author: \href{mailto:yeh@ieee.org}{\nolinkurl{yeh@ieee.org}}}

\begin{document}
\maketitle

\hypertarget{abstract}{%
\section*{Abstract}\label{abstract}}
\addcontentsline{toc}{section}{Abstract}

Multi-pass underfill dispensing couples cycle time with process reliability because the sequence determines the inter-pass flow times. This paper formulates a tri-objective multiset-permutation routing model that minimises nominal cycle time and additive duration variance while maximising reliability under hard- and soft-dwell policies; reliability combines pass, schedule-dependent flow and travel terms. Under independent durations and deterministic dwell thresholds, the variance measure equals the cycle-time variance under soft dwell and bounds it from above under hard dwell. A multiset-permutation binary-addition-tree algorithm enumerates the feasible routes without duplication and prunes prefixes by dominance. Six instances, including three BeagleBone Black benchmarks, yield reference Pareto fronts, the largest search space of \(1.47 \times 10^{12}\) routes in 38 minutes on 20 CPU cores. Among simulated candidates, the nominally cost-selected routes incur maximum estimated excess costs of 0.0014 \$ per board at the modelled duration variability and 0.0035 \$ with it doubled. Round-robin dispatch is within 3.7 \% of the reference optimum at the tested costs; NSGA-II and simplified swarm optimisation show at most 0.53 \% observed cost regret while recovering about 30 \% of the largest hard-dwell front. The framework quantifies sequencing trade-offs and heuristic decision regret; physical reliability prediction requires inspection-based calibration.

\textbf{Keywords:} binary-addition-tree algorithm; multiset permutation; Pareto front; time-reliability optimization; traveling salesman problem; underfill dispensing

\hypertarget{introduction}{%
\section{1. Introduction}\label{introduction}}

Flip-chip and chip-scale packages are attached to printed circuit boards through arrays of solder bumps, and the mismatch in thermal expansion between silicon and the organic substrate concentrates cyclic strain in those bumps. Capillary underfill, an epoxy dispensed at the edge of the die and drawn beneath it by capillary action, redistributes that strain and is the principal means by which the fatigue life of the interconnect is extended {[}1,2{]}. The reliability benefit is, however, conditional on the fillet being formed correctly: trapped air produces voids, and voids under the die act as crack initiation sites during thermal cycling and as moisture reservoirs during reflow {[}2,3{]}. Underfill quality is therefore a reliability question in the strict sense, and the manufacturing step that creates it deserves the same treatment as the components it protects.

For small dies a single pass along one edge, the ``I'' pattern, delivers the whole amount. For larger dies, higher volumes or thinner gaps, the amount is split into several passes: the material from one pass must be allowed to flow under the die before the next pass is applied, otherwise the second amount seals the edge before the first has advanced and a void is trapped {[}4,5{]}. Boards carry several ICs of different sizes, each needing a different number of passes with different amounts and different flow times, and a single dispensing head serves all of them. The head cannot usefully wait at a die while its fillet flows, so it moves to another die and returns later. The sequence in which it does so is not fixed by the process recipe; it is a decision, and it has three consequences. It sets the cycle time of the board, and hence the throughput of the line. It sets the variability of that cycle time, through the travel and dispensing times it accumulates. And, because the gap between consecutive passes on the same die is whatever the sequence leaves, it sets the probability that each fillet is void-free. The same sequencing question arises in die-level underfill, where a flip-chip die is underfilled on its substrate, and in board-level underfill, where ball-grid-array and chip-scale packages are underfilled on the printed circuit board; the model abstracts both, the benchmarks of Section 6 are board-level, and ``IC'' denotes the underfilled device, die or package.

Viewed as a routing problem, the head's tour is a traveling salesman problem (TSP) with two features that the classical and stochastic TSP (STSP) literature does not cover together. First, each ``city'' is visited several times, in a prescribed order and with prescribed amounts, so the feasible tours are permutations of a multiset with a fixed within-die order rather than permutations of distinct nodes. Second, and more fundamentally, the reliability of the tour is not a property of the arcs it uses. In the time-reliability STSP (TR-STSP) {[}6{]}, each arc carries a constant operational probability and the tour reliability is their product; time and reliability interact only through which arcs are chosen. In the dispensing problem, the reliability-carrying quantity is the gap between two visits to the same die, which is a function of the entire schedule the tour induces. The same set of passes, executed in two different orders, has the same travel-time budget to within a few seconds but very different void probabilities. Time and reliability are thus coupled through the schedule, and neither can be optimised in isolation: the fastest tour returns to each die immediately and is essentially guaranteed to void, while a tour that respects every nominal flow time by waiting can more than double the cycle time.

This paper addresses that coupling. Its contributions fall under three headings.

\begin{enumerate}
\def\labelenumi{\arabic{enumi}.}
\tightlist
\item
  \emph{Reliability modelling.} The time-reliability underfill dispensing sequencing problem (TR-UDSP) is defined as a multiset-permutation TSP with precedence. On the nominal schedule, its route reliability, the probability that dispensing succeeds, is the product of a pass reliability depending on the amount, a flow reliability that is a non-decreasing function of the schedule-induced inter-pass gap (parametric here, to be calibrated on void-inspection data), and a travel reliability with a small hazard rate. Hard-dwell and soft-dwell policies are formulated, and the model's boundaries are tested by simulated stochastic execution, by a decomposition of reliability into its route-independent and sequence-dependent parts, and by partial-dwell policies.
\item
  \emph{Exact optimisation.} The permutation BAT {[}6{]} is extended to a multiset-permutation BAT whose selection decoding makes digit prefixes fix route prefixes, so that only feasible routes are generated and prefixes can be pruned by dominance. Completeness, exactness of the objective front under pruning and the work per route are proved, and the method is compared with full enumeration and with best-first branch and bound under the same evaluator, bound and archive.
\item
  \emph{Decision assurance.} Six instances, three derived from the open-hardware BeagleBone Black layout, are solved to reference Pareto fronts of up to \(1.47 \times 10^{12}\) routes. An expected cost per board, or a reliability requirement, turns the front into a selection; against the reference fronts, practical dispatch rules and two metaheuristics (NSGA-II and a multi-objective simplified swarm optimisation, 30 seeds, two budgets) are measured by the part of the front they recover and by the decision regret they incur.
\end{enumerate}

The remainder of the paper is organised as follows. Section 2 gives the notation and assumptions. Section 3 reviews the related work on underfill, TSP variants with multiple visits and precedence, multi-objective optimisation and the BAT family. Section 4 defines the TR-UDSP, its objectives and the cost model. Section 5 presents the multiset-permutation BAT with its correctness and complexity. Section 6 reports the computational results, including the validation of the nominal schedule by simulation and the comparison with dispatch rules, and Section 7 concludes.

\hypertarget{preliminaries}{%
\section{2. Preliminaries}\label{preliminaries}}

\hypertarget{acronyms-notation-and-nomenclature}{%
\subsection{2.1. Acronyms, notation and nomenclature}\label{acronyms-notation-and-nomenclature}}

Tables 1 and 2 list the acronyms and the symbols used throughout; the terms below are used in their stated sense.

\textbf{Table 1.} Acronyms used in the paper.

\begin{longtable}[]{@{}
  >{\raggedright\arraybackslash}p{(\columnwidth - 2\tabcolsep) * \real{0.2189}}
  >{\raggedright\arraybackslash}p{(\columnwidth - 2\tabcolsep) * \real{0.7811}}@{}}
\toprule\noalign{}
\begin{minipage}[b]{\linewidth}\raggedright
Acronym
\end{minipage} & \begin{minipage}[b]{\linewidth}\raggedright
Meaning
\end{minipage} \\
\midrule\noalign{}
\endhead
\bottomrule\noalign{}
\endlastfoot
BAT & binary-addition-tree algorithm \\
CP-SAT & constraint-programming solver of Google OR-Tools {[}7{]} \\
IC & integrated circuit (the underfilled die or package) \\
MP-BAT & multiset-permutation BAT (this paper) \\
MSSO & multi-objective simplified swarm optimisation \\
NSGA-II & non-dominated sorting genetic algorithm II \\
PCB & printed circuit board \\
PLS & Pareto local search \\
TSP / STSP & traveling salesman problem / stochastic TSP \\
TR-STSP & time-reliability STSP {[}6{]} \\
TR-UDSP & time-reliability underfill dispensing sequencing problem \\
\end{longtable}

\textbf{Table 2.} Notation.

\begin{longtable}[]{@{}
  >{\raggedright\arraybackslash}p{(\columnwidth - 2\tabcolsep) * \real{0.2576}}
  >{\raggedright\arraybackslash}p{(\columnwidth - 2\tabcolsep) * \real{0.7424}}@{}}
\toprule\noalign{}
\begin{minipage}[b]{\linewidth}\raggedright
Symbol
\end{minipage} & \begin{minipage}[b]{\linewidth}\raggedright
Meaning
\end{minipage} \\
\midrule\noalign{}
\endhead
\bottomrule\noalign{}
\endlastfoot
\(I = \{ 1,\ldots,m\}\) & set of ICs on the board \\
\(n_{i}\) & number of dispense passes required by IC \(i\) \\
\(q_{i,k}\) & amount dispensed in pass \(k\) of IC \(i\), \(k = 1,\ldots,n_{i}\) \\
\(s_{i},\ e_{i},\ L_{i}\) & start point, end point and length of the dispense line of IC \(i\) \\
\(h\) & home (depot) position of the dispensing head \\
\(v,\ \rho\) & head travel speed; dispense rate \\
\(\tau_{i}\) & nominal capillary flow (dwell) time of IC \(i\) \\
\(N\) & set of pass nodes \((i,k)\); \(|N| = \sum_{i}^{}n_{i}\) \\
\(\iota(u)\) & IC of node \(u\) \\
\(p_{u}\) & processing time of node \(u\): \(p_{i,k} = q_{i,k}/\rho + L_{i}/v\) \\
\(t_{u,w}\) & travel time from the end point of node \(u\) to the start point of node \(w\) \\
\(\sigma( \cdot )\) & standard deviation of a random travel or processing time \\
\(X = \left( x_{1},\ldots,x_{|N|} \right)\) & a route: an ordering of \(N\) executed from \(h\) and returning to \(h\) \\
\(\Omega\) & set of feasible routes \\
\(S(u),\ C(u)\) & start and completion time of node \(u\) under route \(X\) \\
\(g_{i,k}\) & gap \(S(i,k + 1) - C(i,k)\) between consecutive passes on IC \(i\) induced by \(X\) \\
\(W(X)\) & total waiting time of the head under \(X\) (hard-dwell variant) \\
\(r_{i,k}^{P},\ r_{i,k}^{F}(g),\ r_{u,w}^{T}\) & pass, flow and travel reliabilities \\
\(r_{0},\ \kappa,\ \lambda\) & flow reliability at \(g = \tau_{i}\); \(\kappa = - \ln\left( 1 - r_{0} \right)\); travel hazard rate \\
\(T(X),\ \sigma^{*}(X),\ R(X),\ R^{*}(X)\) & nominal cycle time, additive duration variance, route reliability and its logarithm \\
\(K_{P},\ R_{F}(X),\ R_{T}(X)\) & pass, flow and travel factors of \(R(X)\); \(K_{P}\) is route-independent \\
\(F(X),\ \Phi\) & tri-objective value of \(X\); Pareto front \\
\(c_{t},\ c_{f},\ c_{w},\ \mathcal{C}(X)\) & line cost per second, failure cost, extra idle cost per second; expected cost per board \\
\(Z = \left( z_{1},\ldots,z_{|N|} \right),\ a_{j},\ A_{j}\) & MP-BAT combination vector; radix and available-IC list at position \(j\) \\
\end{longtable}

\textbf{Pass.} One execution of the dispense line of an IC, delivering a prescribed amount. The passes of an IC are ordered and must be executed in that order.

\textbf{Dispense line.} The path the needle follows while dispensing; in the I pattern it is a segment along one edge of the die, entered at \(s_{i}\) and left at \(e_{i}\).

\textbf{Route.} A sequence in which the head executes all passes of all ICs, starting and ending at \(h\). A route is feasible if the passes of every IC appear in their prescribed order.

\textbf{Schedule.} The start and completion times of every node, obtained from a route by a forward pass (Section 4.1).

\textbf{Gap.} The time between the completion of pass \(k\) and the start of pass \(k + 1\) on the same IC, as induced by the schedule. Gaps carry the schedule-dependent part of the reliability.

\textbf{Hard-dwell / soft-dwell variant.} In the hard-dwell variant the head waits whenever a gap would be shorter than \(\tau_{i}\); in the soft-dwell variant it never waits and a short gap lowers the flow reliability.

\textbf{Route reliability.} The probability that all passes deliver acceptable amounts, all fillets flow without voiding before the next pass is applied to them, and all moves are clean.

\textbf{Pareto-optimal route.} A route in \(\Omega\) such that no other feasible route is at least as good in all three objectives and strictly better in one.

\hypertarget{assumptions}{%
\subsection{2.2. Assumptions}\label{assumptions}}

\begin{enumerate}
\def\labelenumi{\arabic{enumi}.}
\tightlist
\item
  The board layout, the number of passes per IC, their order and their amounts are given by the process recipe and are not decision variables; only the interleaving of passes across ICs is decided.
\item
  One dispensing head, already loaded with underfill, serves the whole board; no pick-up or tool change is modelled. Each pass is executed in one uninterrupted stroke from \(s_{i}\) to \(e_{i}\).
\item
  Travel between nodes is a straight-line move at constant speed \(v\); its duration is random with known mean \(t_{u,w}\) and standard deviation \(\sigma\left( t_{u,w} \right)\). Processing times are random with known mean \(p_{u}\) and standard deviation \(\sigma\left( p_{u} \right)\). All travel and processing times are mutually independent. The optimisation evaluates every route on its nominal schedule, built from the mean durations, as in {[}6{]}; Section 6.5 tests this against simulated stochastic execution.
\item
  The schedule is non-delay except for the waiting imposed by the dwell rule: a pass starts as soon as the head reaches its start point and, in the hard-dwell variant, the nominal flow time of that IC has elapsed. Intentional waiting beyond that rule is not a decision; the optimisation is therefore restricted to these two waiting policies, and Section 6.6 examines partial-dwell policies.
\item
  Conditionally on the realised schedule, the three failure mechanisms, an out-of-tolerance amount, an insufficiently flowed fillet, and a dirty move, are independent of each other and across passes and moves, so the reliability of a route executed on a given schedule is the product of the individual reliabilities; on the nominal schedule this product is \(R(X)\).
\item
  The flow reliability \(r_{i,k}^{F}\) depends only on the gap \(g_{i,k}\) and on the IC (through \(\tau_{i}\)); it is a known non-decreasing function. In practice it is calibrated from void-inspection data for the die and material in use; the parametric form of Section 4.2 is used here.
\item
  The pass reliability \(r_{i,k}^{P}\) is a known non-increasing function of the amount \(q_{i,k}\), and the travel reliability \(r_{u,w}^{T}\) is a known non-increasing function of the move duration.
\item
  The cost parameters \(c_{t},c_{f},c_{w}\) are known and constant over the planning horizon; the cost of a failed board does not depend on which fillet failed.
\item
  All boards of a production run use the same route, computed offline; adaptive re-sequencing from in-line inspection is left to future work.
\item
  After its last pass, each fillet completes its flow while the board is transferred to and held before curing, outside the dispensing cell; the route reliability therefore contains the \(n_{i} - 1\) inter-pass flows of each IC and not the final one.
\end{enumerate}

\hypertarget{related-work}{%
\section{3. Related Work}\label{related-work}}

Underfill materials, processes and their reliability consequences are reviewed in {[}1,2{]}; the flow of the fillet under the die is a capillary process whose duration grows with die size and gap height and depends on material viscosity and temperature {[}2,3{]}. Industry guidance on automated dispensing describes the multi-pass approach for larger dies and the need to allow flow between passes to avoid entrapped voids {[}4{]}; flow simulation of ball-grid-array underfill confirms that the dispensed weight and the interval between dispenses govern the wet-out and the edge effects {[}5{]}, and the dispensing pattern and method are known to change void formation in capillary underfill {[}8{]}. This literature treats the flow time of one die in isolation; the interaction between several dies served by one head, which is the subject of this paper, is not addressed.

Reliability of a manufacturing system as a whole, the probability that a network of multistate stations with defect rates and rework meets a demand, has been evaluated by Lin and co-workers {[}9,10{]}, and extended by simulation to networks with shared buffers together with a sensitivity analysis of which station's capacity limits the system {[}11{]}; the present paper looks inside one such station and asks for the reliability of the route its tool follows. The nearest analogue is the evaluation of a routing scheme under a time requirement in multistate computer networks {[}12{]}, where the route is given and its reliability is computed; here the route is the decision and its reliability is created by the schedule it induces. Cost is the other attribute usually set against reliability in such networks; Niu et al.~{[}13{]} compute the reliability of a multistate flow network under a cost constraint by enumerating its \((d,b)\)-minimal paths, folding the cost into the reliability definition, whereas the cost model of Section 4.3 is kept outside the objectives and applied to the reference front, so that the same front serves any cost ratio.

The TSP and STSP literature relevant to the TR-STSP is surveyed in {[}6{]} and is not repeated here. Three variants are closer to the present problem. In the sequential ordering problem {[}14{]}, precedence constraints restrict the order of visits; here the precedence is within each IC only, which gives the feasible set the structure of a multiset permutation. In scheduling with minimum time lags {[}15{]}, a job may not start before a fixed delay after its predecessor; the hard-dwell variant is of this type, with a sequence-dependent setup (the travel) and a single machine (the head). Generalised and multiple-visit TSPs are surveyed in {[}16{]}. None of these attaches a reliability to the time lag itself, which is what the soft-dwell variant does.

Pareto optimality, exact and approximate fronts, and the distinction between front-generating methods such as NSGA-II {[}17{]} or multi-objective simplified swarm optimisation {[}18,19{]} and single-solution scalarisations are described in {[}6{]}. Simplified swarm optimisation has been used for bi-objective reliability--redundancy allocation {[}19{]} and for the general active redundancy allocation problem {[}20{]}, and has tuned the hyper-parameters of convolutional neural networks {[}21{]}; the first two are the reliability-design counterparts of the present sequencing problem: a discrete decision, a reliability objective in product form, and a cost that is traded against it. The BAT was introduced for binary-state network reliability {[}22{]} and accelerated by early termination in the quick BAT {[}23{]}. The permutation BAT of {[}6{]} replaces the binary digits by mixed-radix digits to enumerate permutations of distinct elements with an inversion-table decoding. The present work generalises that counter to permutations of a multiset with fixed within-class order and replaces the inversion-table decoding by a selection decoding, under which the first \(p\) digits fix the first \(p\) stations of the route whatever the remaining digits are. That property, which the inversion-table decoding lacks because a later digit shifts earlier elements, is what allows prefixes to be pruned by dominance. Gray-code orders of multiset permutations by interchanges {[}24{]} are the other classical enumeration order and are discussed in Section 5.2. The rank decoding used here coincides, when all multiplicities are one, with the factorial number system of {[}25{]}.

\hypertarget{the-time-reliability-underfill-dispensing-sequencing-problem}{%
\section{4. The Time-Reliability Underfill Dispensing Sequencing Problem}\label{the-time-reliability-underfill-dispensing-sequencing-problem}}

\hypertarget{process-setting-feasible-routes-and-induced-schedule}{%
\subsection{4.1. Process setting, feasible routes and induced schedule}\label{process-setting-feasible-routes-and-induced-schedule}}

A board carries \(m\) ICs. IC \(i\) requires \(n_{i}\) passes with prescribed amounts \(q_{i,1},\ldots,q_{i,n_{i}}\) along its dispense line from \(s_{i}\) to \(e_{i}\), and a nominal flow time \(\tau_{i}\) between consecutive passes. A single head starts at \(h\), executes all \(\sum_{i}^{}n_{i}\) passes, and returns to \(h\). The processing time of node \((i,k)\) is \(p_{i,k} = q_{i,k}/\rho + L_{i}/v\) and the travel time between nodes is \(t_{u,w} = \parallel e_{\iota(u)} - s_{\iota(w)} \parallel /v\). Figure 1 shows the two board layouts used in Section 6.

\includegraphics[width=6in,height=2.52in]{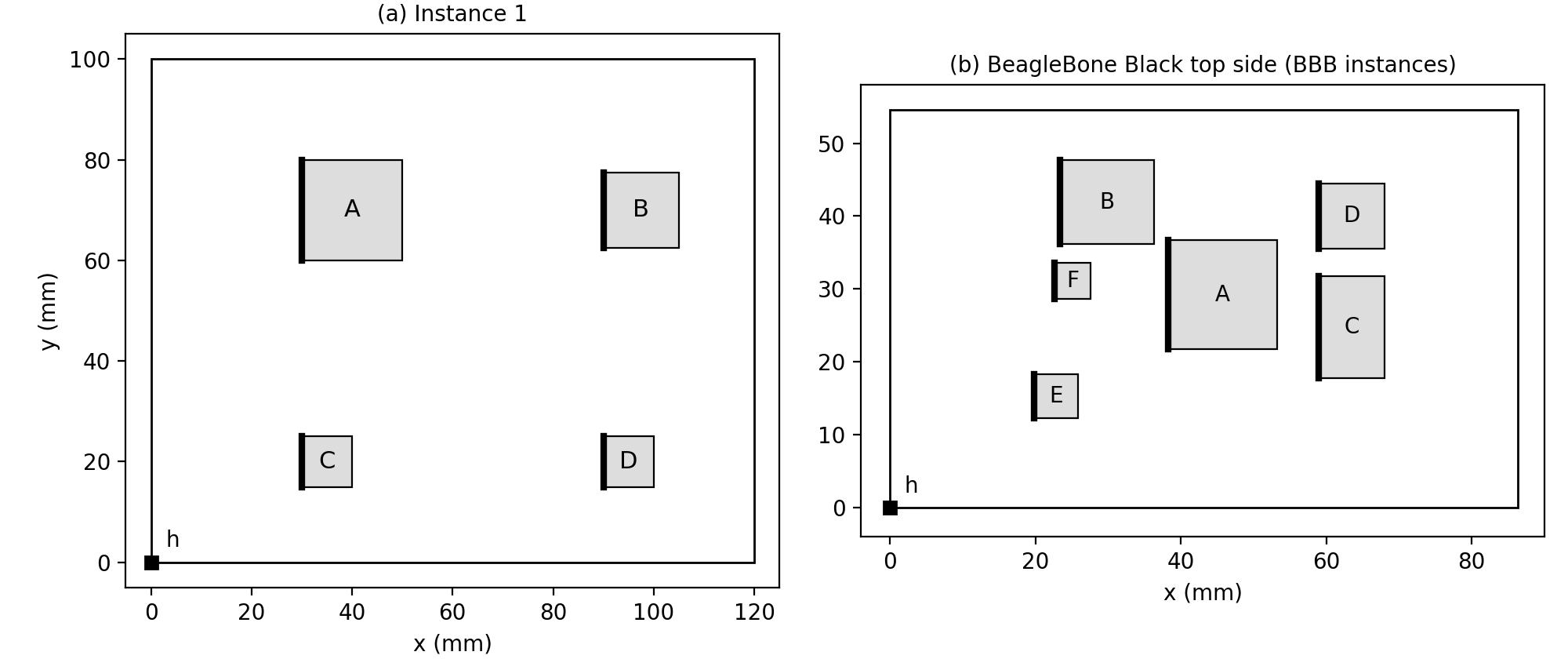}

Figure 1. Board layouts. (a) Synthetic Instance 1 (120 × 100 mm, four ICs). (b) Top side of the BeagleBone Black Rev C (86.4 × 54.6 mm, six ICs). Thick left edges mark the I-shaped dispense lines; \(h\) is the home position.

A route \(X\) is feasible if, for every IC \(i\), its passes appear in the order \((i,1),(i,2),\ldots,\left( i,n_{i} \right)\). The set \(\Omega\) of feasible routes is the set of permutations of the multiset \(\{ 1^{n_{1}},\ldots,m^{n_{m}}\}\), whose cardinality is

\begin{equation}|\Omega| = \frac{\left( \sum_{i}^{}n_{i} \right)!}{\prod_{i}^{}n_{i}!}.\end{equation} For Instance 1 (\(n = 5,4,3,3\)), \(|\Omega| = 15!/(5!\, 4!\, 3!\, 3!) = 12,612,600\), whereas the \(15!\) orderings of the pass nodes number \(1.3 \times 10^{12}\). Given \(X\), the schedule is computed by one forward pass,

\begin{equation}S\left( x_{1} \right) = t_{h,x_{1}},\quad\quad S\left( x_{j} \right) = C\left( x_{j - 1} \right) + t_{x_{j - 1},x_{j}},\quad\quad C\left( x_{j} \right) = S\left( x_{j} \right) + p_{x_{j}},\end{equation} and the induced gap of IC \(i\) between passes \(k\) and \(k + 1\) is \(g_{i,k} = S(i,k + 1) - C(i,k)\). In the hard-dwell variant the head waits, if necessary, so that \(g_{i,k} \geq \tau_{i}\), and the waiting time \(W(X)\) is added to the schedule; in the soft-dwell variant the head never waits and short gaps are penalised through the reliability model.

\hypertarget{objectives-time-and-schedule-dependent-reliability}{%
\subsection{4.2. Objectives: time and schedule-dependent reliability}\label{objectives-time-and-schedule-dependent-reliability}}

Following {[}6{]}, the first two objectives are the nominal cycle time, computed from the mean durations, and the additive duration variance \(\sigma^{*}\), the sum of the variances of the durations on the route; both are additive along the route:

\begin{equation}\min\ T(X) = t_{h,x_{1}} + \sum_{j = 1}^{|N| - 1}t_{x_{j},x_{j + 1}} + t_{x_{|N|},h} + \sum_{(i,k) \in N}^{}p_{i,k} + W(X),\end{equation} \begin{equation}\min\ \sigma^{*}(X) = \sum_{(u,w) \in X}^{}\sigma\left( t_{u,w} \right)^{2} + \sum_{(i,k) \in N}^{}\sigma\left( p_{i,k} \right)^{2}.\end{equation} The second sum in Equation (4) does not depend on \(X\) and is kept for reporting only. Variability is an objective in its own right, not a by-product of the mean, for a reason that belongs to the line rather than to the board: a dispensing cell runs to a takt, and the buffer before the next station is sized by how much the cycle time varies, not by its mean, so a line needs a measure of variability that it can bound. Two routes with the same cycle time are not equivalent if one of them spends its time in a few long moves and the other in many short ones; the first accumulates more duration variance. The pass part of \(\sigma^{*}\) is the same for every route, so what \(\sigma^{*}\) discriminates is the travel. What \(\sigma^{*}\) means for the realised cycle time depends on the dwell policy, and Proposition 1 states it exactly. In the soft-dwell variant the realised cycle time is the sum of the route's durations, so \(T\) is its mean and \(\sigma^{*}\) its variance. In the hard-dwell variant each wait is the maximum of two random quantities, the arrival of the head and the end of the dwell; the realised cycle time then has a mean at least \(T\) and a variance at most \(\sigma^{*}\), because the waits absorb part of the variability accumulated before them. Minimising the bound does not minimise the realised variance, and the gap between the two is measured in Section 6.5. The squared form is used so that the objective is additive along a prefix, which is what the bound of Section 5.2 requires; the tables report its square root, which is the quantity the line engineer reads.

\textbf{Proposition 1.} Fix a route \(X\) and let \(D = \left( D_{1},\ldots,D_{K} \right)\) be its travel and processing durations, including the trips from and to \(h\), mutually independent with finite variances and with means equal to the nominal durations. Let \(\widetilde{T}(X)\) be the realised cycle time of \(X\) under the hard-dwell rule with the deterministic thresholds \(\tau_{i}\). Then \(T(X) \leq \mathbb{E}\left\lbrack \widetilde{T}(X) \right\rbrack\) and \(\mathrm{Var}\left\lbrack \widetilde{T}(X) \right\rbrack \leq \sigma^{*}(X)\). Under the soft-dwell rule both hold with equality.

\emph{Proof.} Under the hard-dwell rule the start of each node is the later of the arrival of the head and the end of the dwell of that IC, and every completion is a start plus a processing time, so the recursion of Section 4.1 consists of sums and maxima only. Unrolling it, \(\widetilde{T}(X) = f_{X}(D) = \max_{\pi}\left( a_{\pi} + \sum_{j \in \pi}D_{j} \right)\), the longest path in the acyclic graph whose nodes are the start and the completion of every pass and whose arcs are the passes (start to completion, weighted by the processing duration), the moves of the route (completion to the next start, weighted by the travel duration, with the trips from and to \(h\)) and the dwells (completion of pass \(k\) of IC \(i\) to the start of its pass \(k + 1\), weighted by the constant \(\tau_{i}\)); \(\pi\) ranges over the finite set of paths from the departure from \(h\) to the return to \(h\), \(a_{\pi}\) is the sum of the thresholds on \(\pi\), and every path uses each duration at most once because the positions along a path increase. First, \(f_{X}\) is a maximum of affine functions of \(D\) and therefore convex, so \(T(X) = f_{X}\left( \mathbb{E}\lbrack D\rbrack \right) \leq \mathbb{E}\left\lbrack f_{X}(D) \right\rbrack\) by Jensen's inequality. Second, replacing one duration \(D_{j}\) by an independent copy \(D_{j}^{\prime}\) changes every affine term by at most \(\left| D_{j} - D_{j}^{\prime} \right|\) and hence the maximum by at most that amount; the Efron--Stein inequality {[}26,27{]} then gives \(\mathrm{Var}\left\lbrack f_{X}(D) \right\rbrack \leq \frac{1}{2}\sum_{j}\mathbb{E}\left\lbrack \left( f_{X}(D) - f_{X}\left( D^{(j)} \right) \right)^{2} \right\rbrack \leq \frac{1}{2}\sum_{j}\mathbb{E}\left\lbrack \left( D_{j} - D_{j}^{\prime} \right)^{2} \right\rbrack = \sum_{j}\mathrm{Var}\left( D_{j} \right) = \sigma^{*}(X)\), where \(D^{(j)}\) is \(D\) with the \(j\)-th coordinate replaced. Under the soft-dwell rule no maximum is taken and \(\widetilde{T}(X) = \sum_{j}D_{j}\), whose mean is \(T(X)\) and whose variance is \(\sigma^{*}(X)\). □

The reliability of a route is the probability that it is executed successfully: that every pass is in tolerance, every fillet has flowed before the next pass is applied to it, and every move is clean, up to the moment the board leaves the dispensing cell; the flow of the final pass of each IC completes outside the cell (Assumption 10). It represents modelled dispensing-process success, not a prediction of solder-joint fatigue life or field reliability, to which it is an input. Three independent failure mechanisms are modelled. (a) Pass reliability: each pass may deliver an out-of-tolerance amount because of viscosity drift, partial nozzle clogging or tailing; \(r_{i,k}^{P} = r^{P}\left( q_{i,k} \right)\) is non-increasing in the amount and independent of \(X\). (b) Flow reliability: pass \(k + 1\) on IC \(i\) is applied on a fillet that has had \(g_{i,k}\) seconds to flow; the probability that the flow has progressed far enough to avoid a void is a non-decreasing function of the gap,

\begin{equation}r_{i,k}^{F}\left( g_{i,k} \right) = 1 - \exp\left( - \kappa\,\frac{g_{i,k}}{\tau_{i}} \right),\quad\quad\kappa = - \ln\left( 1 - r_{0} \right),\end{equation} so that \(r^{F} = r_{0}\) at the nominal dwell, \(r^{F} \rightarrow 1\) for longer gaps and \(r^{F} \rightarrow 0\) as the gap vanishes; any other non-decreasing function, for example one fitted to inspection data, can be substituted without changing the algorithm. (c) Travel reliability: a move between two dies is not risk-free. The needle travels loaded, at the dispense height, over other components; the longer the move, the greater the chance that a drop detaches or a string is drawn onto a neighbouring part, that a satellite is shed by the acceleration and deceleration of the head, or that the accumulated positioning error at the far end of a long move leaves the needle off the die edge at the start of the next pass. Each of these is a failure of the board that occurs during the move rather than during a pass, and each becomes more likely the longer the move lasts. The simplest model with these properties is a constant hazard rate \(\lambda\) per second of travel, \(r_{u,w}^{T} = \exp\left( - \lambda t_{u,w} \right)\), which is the arc reliability of the TR-STSP {[}6{]} given a physical meaning: it is a property of the move, not of the schedule, and its logarithm is proportional to travel time, so it enters Equation (7) as a term that is additive along the route and pulls in the same direction as \(T\). With the value of \(\lambda\) used here it changes \(R\) by a few tenths of a percent over a whole route; it is retained because it is the term that distinguishes two routes with identical passes and gaps but different travel, and because a line with a stringing-prone material can raise \(\lambda\) without changing anything else. Under Assumption 5 the route reliability is

\begin{equation}R(X) = \prod_{(i,k) \in N}^{}r_{i,k}^{P} \cdot \prod_{i \in I}^{}{\prod_{k = 1}^{n_{i} - 1}r_{i,k}^{F}}\left( g_{i,k}(X) \right) \cdot \prod_{(u,w) \in X}^{}r_{u,w}^{T},\end{equation} and, as in {[}6{]}, it is maximised in logarithmic form so that it is additive along the route:

\begin{equation}\max\ R^{*}(X) = \sum_{(i,k) \in N}^{}\ln r_{i,k}^{P} + \sum_{i \in I}^{}{\sum_{k = 1}^{n_{i} - 1}\ln}r_{i,k}^{F}\left( g_{i,k}(X) \right) - \lambda\sum_{(u,w) \in X}^{}t_{u,w}.\end{equation} The first sum is a constant of the instance. The second couples \(R^{*}\) to the schedule: shortening the cycle time by visiting the same IC again too early lowers \(g_{i,k}\) and therefore \(R^{*}\), while lengthening the gaps raises \(R^{*}\) at the cost of \(T\). In the hard-dwell variant \(g_{i,k} \geq \tau_{i}\) always holds and the second sum is bounded below by \(\sum_{i}^{}\left( n_{i} - 1 \right)\ln r_{0}\); in the soft-dwell variant it can become arbitrarily negative. Writing \(K_{P} = \prod r_{i,k}^{P}\) for the route-independent pass factor and \(R_{F}(X)\), \(R_{T}(X)\) for the flow and travel products, \(R(X) = K_{P}\, R_{F}(X)\, R_{T}(X)\): \(K_{P}\) bounds \(R\) from above for every route (0.811 on Instance 1), and \(R/K_{P}\) is the part that sequencing controls. \(R\) is evaluated at the nominal gaps. Because \(r^{F}\) is concave in the gap, each flow factor taken alone averages downward under random durations, \(\mathbb{E}\left\lbrack r^{F}(G) \right\rbrack \leq r^{F}\left( \mathbb{E}\lbrack G\rbrack \right)\); this does not fix the sign of the route-level error, because the realised gaps and travel factors share durations and are dependent, and under hard dwell the nominal gaps are not the means of the realised ones. Section 6.5 therefore estimates the product of the realised factors directly.

The TR-UDSP is therefore

\begin{equation}\min\ T(X),\quad\quad\min\ \sigma^{*}(X),\quad\quad\max\ R^{*}(X),\quad\quad X \in \Omega,\end{equation} and its solution is the Pareto front \(\Phi\) of Equation (8). Table 3 contrasts it with the TR-STSP.

\textbf{Table 3.} TR-UDSP versus TR-STSP.

\begin{longtable}[]{@{}
  >{\raggedright\arraybackslash}p{(\columnwidth - 4\tabcolsep) * \real{0.2692}}
  >{\raggedright\arraybackslash}p{(\columnwidth - 4\tabcolsep) * \real{0.2186}}
  >{\raggedright\arraybackslash}p{(\columnwidth - 4\tabcolsep) * \real{0.5122}}@{}}
\toprule\noalign{}
\begin{minipage}[b]{\linewidth}\raggedright
Aspect
\end{minipage} & \begin{minipage}[b]{\linewidth}\raggedright
TR-STSP {[}6{]}
\end{minipage} & \begin{minipage}[b]{\linewidth}\raggedright
TR-UDSP (this paper)
\end{minipage} \\
\midrule\noalign{}
\endhead
\bottomrule\noalign{}
\endlastfoot
Node visits & each city once & IC \(i\) visited \(n_{i}\) times, in fixed order \\
Solution space & \((n - 1)!\) permutations & \(\left( \sum_{i}^{}n_{i} \right)!/\prod_{i}^{}n_{i}!\) multiset permutations \\
Reliability carrier & arc, constant & pass (constant) + inter-pass gap (schedule-dependent) + arc \\
Time--reliability coupling & through arc choice only & through the induced schedule \(g_{i,k}(X)\) \\
Dwell & none & hard lag, or soft (reliability-penalised) \\
\end{longtable}

\hypertarget{economic-interpretation-expected-cost-per-board}{%
\subsection{4.3. Economic interpretation: expected cost per board}\label{economic-interpretation-expected-cost-per-board}}

The front \(\Phi\) is what the optimisation returns; the line engineer must pick one route. Each route is therefore also assigned an expected cost per board,

\begin{equation}\mathcal{C}(X) = c_{t}\, T(X) + c_{f}\,\left( 1 - R(X) \right) + c_{w}\, W(X),\end{equation} where \(c_{t}\) is the line cost per second of cycle time, \(c_{f}\) the expected cost of a board whose underfill fails, and \(c_{w}\) an extra cost per second of head waiting beyond the line cost, which is already charged through \(T\) because \(T\) contains \(W\) (Equation (3)); \(c_{w}\) covers idle-specific costs such as needle-tip drying or pot-life consumed while the head waits, and \(c_{w} = 0\), used in all experiments, recovers the pure line-cost model. With \(c_{w} = 0\), \(\mathcal{C}\) is non-decreasing in \(T\) and non-increasing in \(R\), so every route is weakly dominated in cost by a member of \(\Phi\) and the cheapest member of \(\Phi\) is the cheapest route in \(\Omega\). With \(c_{w} > 0\) this no longer holds, because a route dominated in \((T,\sigma^{*},R)\) may wait less; the selection would then need \(W\) as a fourth criterion. Where a failure cost is hard to state, a reliability requirement selects instead: the fastest member of \(\Phi\) with \(R \geq R_{\min}\), read off the same front. Relative to a current-practice route \(X_{0}\), the saving per board is \(\Delta\mathcal{C}(X) = \mathcal{C}\left( X_{0} \right) - \mathcal{C}(X)\) and the annual saving on a fixed volume of \(V_{b}\) boards per year is

\begin{equation}\Delta\mathcal{C}(X) \times V_{b}.\end{equation} When the line is capacity-constrained, the shorter cycle time also frees \(V_{b}\left( T\left( X_{0} \right) - T(X) \right)\) seconds of line time per year; its value depends on what the line can do with that time and is not added to Equation (10). Equation (9) is not a fourth objective; it is a scalarisation applied after the front is known, so that different cost ratios select different points on the same front without re-running the optimisation, in the decision-support sense of {[}6{]}.

\hypertarget{multiset-permutation-bat}{%
\section{5. Multiset-Permutation BAT}\label{multiset-permutation-bat}}

\hypertarget{combination-vector-and-generation-procedure}{%
\subsection{5.1. Combination vector and generation procedure}\label{combination-vector-and-generation-procedure}}

The permutation BAT {[}6{]} is a mixed-radix BAT with upper bounds \(UB\left( z_{k} \right) = k - 1\), decoded by shifting the \(k\)-th element \(z_{k}\) positions to the left; it presupposes distinct elements that are each visited once. In the TR-UDSP the five passes of IC A are not interchangeable in position. Enumerating the \(15!\) orderings of the pass nodes and discarding those with the wrong within-IC order generates about \(10^{5}\) vectors for every feasible route; a fixed-radix vector \(z_{j} \in \{ 0,\ldots,m - 1\}\) generates \(m^{|N|} = 4^{15} \approx 1.07 \times 10^{9}\) vectors of which 99 \% have wrong letter counts. Neither preserves the defining property of a BAT, that the counter visits only the objects wanted, once each.

A feasible route is fully determined by answering, at each position \(j = 1,\ldots,|N|\), the question ``which IC executes its next pass now?''; the pass index is implied. The route is encoded as \(Z = \left( z_{1},\ldots,z_{|N|} \right)\) where \(z_{j}\) is the rank of the chosen IC in the sorted list \(A_{j}\) of ICs that still have at least one pass remaining before position \(j\). The radix of position \(j\) is \(a_{j} = \left| A_{j} \right|\) and \(z_{j} \in \{ 0,\ldots,a_{j} - 1\}\); decoding is one lookup, \(x_{j} = A_{j}\left\lbrack z_{j} \right\rbrack\), followed by decrementing that IC's remaining count. The digit \(z_{j}\) is a zero-based rank in the current available-IC list, not a fixed identifier of the IC: once an IC has executed all its passes it leaves \(A_{j}\), and the ranks of the ICs after it shift down. The radix therefore starts at \(m\) and shrinks as ICs are exhausted, and it depends on the prefix \(z_{1}\cdots z_{j - 1}\) rather than on \(j\) alone; this is the only structural difference from the mixed-radix counters of the BAT family. The decoding is a different bijection from the inversion-table decoding of {[}6{]}, not a re-indexing of it: under shift decoding the digit \(z_{k}\) moves the \(k\)-th element past earlier ones, so even the first station of the route depends on every digit, whereas here \(z_{1}\cdots z_{p}\) determine the first \(p\) stations regardless of the suffix. Section 5.2 rests on this property. Table 4 lists the vectors for two ICs with two passes each, and Table 5 decodes four vectors of Instance 1.

\textbf{Table 4.} All six combination vectors and routes for ICs A and B with \(n_{A} = n_{B} = 2\).

\begin{longtable}[]{@{}
  >{\raggedright\arraybackslash}p{(\columnwidth - 6\tabcolsep) * \real{0.2723}}
  >{\raggedright\arraybackslash}p{(\columnwidth - 6\tabcolsep) * \real{0.1961}}
  >{\raggedright\arraybackslash}p{(\columnwidth - 6\tabcolsep) * \real{0.3486}}
  >{\raggedright\arraybackslash}p{(\columnwidth - 6\tabcolsep) * \real{0.1830}}@{}}
\toprule\noalign{}
\begin{minipage}[b]{\linewidth}\raggedright
\#
\end{minipage} & \begin{minipage}[b]{\linewidth}\raggedright
\(Z\)
\end{minipage} & \begin{minipage}[b]{\linewidth}\raggedright
\(\left( a_{1},a_{2},a_{3},a_{4} \right)\)
\end{minipage} & \begin{minipage}[b]{\linewidth}\raggedright
Route
\end{minipage} \\
\midrule\noalign{}
\endhead
\bottomrule\noalign{}
\endlastfoot
1 & (0,0,0,0) & 2, 2, 1, 1 & A A B B \\
2 & (0,1,0,0) & 2, 2, 2, 1 & A B A B \\
3 & (0,1,1,0) & 2, 2, 2, 1 & A B B A \\
4 & (1,0,0,0) & 2, 2, 2, 1 & B A A B \\
5 & (1,0,1,0) & 2, 2, 2, 1 & B A B A \\
6 & (1,1,0,0) & 2, 2, 1, 1 & B B A A \\
\end{longtable}

\textbf{Table 5.} Four combination vectors of Instance 1 and their routes, generated by the encoder used in the experiments and checked by decoding.

\begin{longtable}[]{@{}
  >{\raggedright\arraybackslash}p{(\columnwidth - 4\tabcolsep) * \real{0.3723}}
  >{\raggedright\arraybackslash}p{(\columnwidth - 4\tabcolsep) * \real{0.2364}}
  >{\raggedright\arraybackslash}p{(\columnwidth - 4\tabcolsep) * \real{0.3913}}@{}}
\toprule\noalign{}
\begin{minipage}[b]{\linewidth}\raggedright
\(Z\)
\end{minipage} & \begin{minipage}[b]{\linewidth}\raggedright
Route
\end{minipage} & \begin{minipage}[b]{\linewidth}\raggedright
Remark
\end{minipage} \\
\midrule\noalign{}
\endhead
\bottomrule\noalign{}
\endlastfoot
(0,0,0,0,0,0,0,0,0,0,0,0,0,0,0) & AAAAABBBBCCCDDD & first vector: IC-by-IC order \\
(0,0,0,0,0,0,0,0,0,0,0,1,0,0,0) & AAAAABBBBCCDCDD & second vector: first carry, at position 12 \\
(0,1,0,1,3,2,0,1,3,2,0,1,2,1,0) & ABABDCABDCABDCA & a min-\(T\) member of the hard-dwell front \\
(3,3,3,2,2,2,1,1,1,1,0,0,0,0,0) & DDDCCCBBBBAAAAA & last vector: reverse IC-by-IC order \\
\end{longtable}

Algorithm 1 generates every vector once, in lexicographic order, starting from the all-zero vector. It is the BAT carry rule with one addition, line 3: because the radices of positions after the changed digit depend on what remains, the zeroed suffix is rebuilt rather than merely reset. When the bound prunes at depth \(p\), the carry starts at position \(p\) instead of \(|N|\) (lines 5 and 7), which skips exactly the vectors that share the digit prefix \(z_{1}\cdots z_{p}\); positions after \(p\) were never decoded and are rebuilt from the new prefix. Algorithm 1 also embeds the front update and the prefix bound described below; when the bound is inactive it reduces to plain enumeration.

\begin{center}\rule{0.5\linewidth}{0.5pt}\end{center}

\textbf{Algorithm 1.} MP-BAT: multiset-permutation BAT with prefix dominance bound

\textbf{Input:} \(n_{i}\), \(q_{i,k}\), \(\tau_{i}\), coordinates; variant (hard/soft); seed front \(\Phi_{0}\) (possibly empty). \textbf{Output:} Pareto front \(\Phi\).

\begin{enumerate}
\def\labelenumi{\arabic{enumi}.}
\tightlist
\item
  \textbf{initialise} \(\Phi \leftarrow \Phi_{0}\) sorted by \(T\); \({rem}_{i} \leftarrow n_{i}\); \(Z \leftarrow (0,\ldots,0)\); \({state}_{0} \leftarrow (t = 0,\) \(V = \sum\sigma( p_{u} )^{2},\) \(LR = \sum\ln r_{u}^{P})\); \(j \leftarrow 1\).
\item
  pruned \(\leftarrow\) false; \textbf{for} \(p = j,\ldots,|N|\) \textbf{do}
\item
  \(\quad A_{p} \leftarrow\) sorted list of \(i\) with \({rem}_{i} > 0\); \(a_{p} \leftarrow \left| A_{p} \right|\); \(x_{p} \leftarrow A_{p}\left\lbrack z_{p} \right\rbrack\); \({rem}_{x_{p}} \leftarrow {rem}_{x_{p}} - 1\);
\item
  \(\quad{state}_{p} \leftarrow Extend\left( {state}_{p - 1},x_{p} \right)\) (clock, variance, log-reliability, last completion; Section 5.2);
\item
  \(\quad\) \textbf{if} \(p < |N|\) \textbf{and} \(Bound\left( {state}_{p} \right)\) is dominated by \(\Phi\) \textbf{then} \(q \leftarrow p\); mark pruned; \textbf{break} (the block of vectors with prefix \(z_{1}\cdots z_{p}\) is skipped).
\item
  \textbf{if} not pruned \textbf{then} \(F(X) \leftarrow {state}_{|N|}\) plus return trip; insert \(F(X)\) into \(\Phi\), removing members it dominates; \(q \leftarrow |N|\).
\item
  \(j \leftarrow q\); \textbf{while} \(j \geq 1\) \textbf{and} \(z_{j} = a_{j} - 1\) \textbf{do} \({rem}_{x_{j}} \leftarrow {rem}_{x_{j}} + 1\); \(j \leftarrow j - 1\) (carry).
\item
  \textbf{if} \(j = 0\) \textbf{then return} \(\Phi\) \textbf{else} \({rem}_{x_{j}} \leftarrow {rem}_{x_{j}} + 1\); \(z_{j} \leftarrow z_{j} + 1\); \(z_{j + 1..|N|} \leftarrow 0\); \textbf{go to} line 2.
\end{enumerate}

\begin{center}\rule{0.5\linewidth}{0.5pt}\end{center}

\hypertarget{incremental-evaluation-prefix-dominance-bound-warm-start-and-parallel-decomposition}{%
\subsection{5.2. Incremental evaluation, prefix dominance bound, warm start and parallel decomposition}\label{incremental-evaluation-prefix-dominance-bound-warm-start-and-parallel-decomposition}}

Consecutive vectors share a prefix, so the state after positions \(1,\ldots,p - 1\) (clock \(t\), accumulated variance \(V\), log-reliability \(LR\), last completion of each IC, remaining counts) is reused and only positions from the changed digit onward are re-evaluated; Section 5.3 gives the resulting work per route. The initial state already contains the route-independent pass terms of \(\sigma^{*}\) and \(R^{*}\) (Algorithm 1, line 1). Let \(t_{\min}\) be the shortest move of the instance over every pair of ICs and over the trips from and to \(h\); the return to \(h\) is one of the \(|N| - p + 1\) moves still to be made, so \(t_{\min}\) must cover it. Along a prefix, \(t\) and \(V\) are non-decreasing and \(LR\) is non-increasing, so at depth \(p\) every completion of the prefix satisfies

\begin{equation}T \geq T_{lb} = t_{p} + \sum_{u\,\text{remaining}}p_{u} + t_{\min}\,\left( |N| - p + 1 \right),\quad\quad\sigma^{*} \geq V_{p},\quad\quad R^{*} \leq LR_{p},\end{equation} and in the hard-dwell variant \(T_{lb}\) is raised to \(\max_{i:\,\ell_{i} > 0}\left\lbrack \hat{C}_{i} + \hat{\ell}_{i}\tau_{i} + \sum_{k\,\text{remaining}}p_{i,k} + t_{i,h} \right\rbrack\), where \(\ell_{i}\) is the number of passes of IC \(i\) still to execute, \(t_{i,h}\) the travel time from IC \(i\) to \(h\), \(\hat{C}_{i}\) the completion time of the last executed pass of IC \(i\) and \(\hat{\ell}_{i} = \ell_{i}\) if IC \(i\) has started, and \(\hat{C}_{i} = t_{p}\) and \(\hat{\ell}_{i} = \ell_{i} - 1\) otherwise. If some member of \(\Phi\) weakly dominates \(\left( T_{lb},V_{p},LR_{p} \right)\), no completion of the prefix can enter \(\Phi\) with a new objective point, and the block of vectors with that prefix is skipped by starting the carry at position \(p\): the skip mechanism of the quick BAT {[}23{]}, driven by dominance rather than by feasibility. Two facts make the skip legitimate, and both are needed: the digit prefix fixes the route prefix (Section 5.1), so every vector in the skipped block starts with the same \(p\) stations; and the three bounds are monotone along the prefix, so what holds at depth \(p\) holds for every completion. The pruned block is contiguous in the counter because the enumeration is lexicographic. The archive is kept sorted by \(T\) so that a dominance test stops at the first member with \(T > T_{lb}\).

\emph{Floating-point implementation.} Propositions 2 and 3 (Section 5.3) establish complete enumeration and front preservation in exact arithmetic. The implementation tests dominance with absolute tolerances of \(10^{-9}\) s on \(T\) and \(10^{-12}\) on \(\sigma^{*}\) and \(R^{*}\), in pruning and in archive insertion alike; the tolerances therefore decide which prefixes are pruned and which archive members are deleted, not only how near-equal points are reported, and the relaxed relation is not transitive, so a chain of replacements could in principle drift by more than one tolerance. The evidence that this does not affect the reported fronts is computational: on every instance enumerated in full, the pruned front coincides with the full-enumeration front, which uses the same insertion tolerance but no bound (Section 6.3, Table 12); the serial and parallel programs agree (Section 6.1, control (ii)); and the CP-SAT control agrees at the reported precision (control (iv)). These checks are not interval-arithmetic certificates, which is why the computed fronts are called reference fronts (Section 6.1); a comparison that falls within a tolerance can be recomputed in higher precision, and none of the reported selections depends on one.

Seeding \(\Phi_{0}\) with a heuristic front does not change the result of a completed run but makes the bound active from the first vector; Section 6.3 shows that this decides between failure and progress at \(10^{12}\) routes. Fixing the first \(d\) digits fixes a contiguous block of the counter; blocks are independent sub-counters with their own archives, merged by dominance at the end, which is the same communication-free split used in parallel BAT implementations.

\emph{Remark (order of enumeration).} Lexicographic order is what makes a pruned prefix a contiguous block of the counter. A Gray-code order of multiset permutations by interchanges {[}24{]} changes two positions between consecutive routes, but in this problem an interchange moves every later start time, gap and hard-dwell wait, so the objectives cannot be updated in constant time, and routes sharing a prefix are scattered through the sequence, so no block can be skipped. Best-first branch and bound expands prefixes in increasing \(T_{lb}\) and stores the open set; Section 6.3 compares it with the counter under the same evaluator, bound and archive. A value ordering of \(A_{j}\) by travel time from the last executed IC is compatible with the counter and would make the archive tighten earlier; it is left for future work.

\hypertarget{correctness-and-complexity}{%
\subsection{5.3. Correctness and complexity}\label{correctness-and-complexity}}

\textbf{Proposition 2.} The map from feasible routes to digit sequences with \(z_{j} < a_{j}\) is a bijection, and Algorithm 1 without pruning outputs every feasible route exactly once.

\emph{Proof.} Encoding a route by recording, at each position, the rank of the executed IC among those with passes remaining is a deterministic function; decoding inverts it position by position, so the map is injective and its image is exactly the set of admissible digit sequences, whose number is the multinomial of Equation (1). Lines 7--8 implement the successor function of the mixed-radix counter over admissible sequences: the least significant digit below its maximum is incremented and all lower digits are reset to their minimum, after which lines 2--3 rebuild the suffix consistently with the new prefix. Starting from the minimum sequence and stopping when no digit can be incremented, the counter visits each admissible sequence once. □

\textbf{Proposition 3.} With the bound of Equation (11), Algorithm 1 outputs the same set of Pareto-optimal objective points as without pruning, and for each point a route attaining it.

\emph{Proof.} \(T\) and \(\sigma^{*}\) accumulate non-negative terms and \(R^{*}\) non-positive terms along a prefix, and the route-independent constants are in the initial state, so every completion of a pruned prefix is weakly dominated in all three objectives by the archive member that dominated the bound. The archive replaces a member only by one that dominates it, so at the end every pruned completion is weakly dominated by a member of \(\Phi\): it is either dominated, and not Pareto-optimal, or it attains an objective point that \(\Phi\) already represents. Routes with exactly equal objectives may therefore be pruned; the front is exact as a set of objective points, each with one attaining route, not as a list of all attaining routes. □

Without pruning, Algorithm 1 evaluates each distinct feasible prefix once. Let \(P_{p}\) be the number of distinct feasible prefixes of length \(p\) and \(P = \sum_{p < |N|}P_{p}\); the work is \(O\left( m\,(P + |\Omega|) + |\Omega|\,|\Phi| \right)\), the last term for the archive tests. Every prefix has at least one completion, so \(P \leq |N|\,|\Omega|\); \(P/|\Omega|\) is the number of internal prefixes evaluated per feasible route, and the positions evaluated per route, leaves included, number \(\left( P + |\Omega| - 1 \right)/|\Omega|\), one more than \(P/|\Omega|\) less the root. No constant bound on \(P/|\Omega|\) holds in general: for two ICs with multiplicities \((|N| - 1,1)\) there are \(p + 1\) prefixes of length \(p\) and \(|N|\) routes, so \(P/|\Omega| = (|N| + 1)/2\). Balanced multiplicities keep it small. For the six instances of Section 6, \(P/|\Omega|\) is 1.88--2.06 (Instance 1: 25,959,593 prefixes for 12,612,600 routes, equal to the count of the instrumented solver), against 2.22 for \((5,5,5)\), 2.78 for \((8,8)\), 4.1 for \((12,1,1,1)\), 8.0 for \((14,1)\) and 24.5 for \((47,1)\). With pruning the number of evaluated prefixes depends on the instance and the seed and is reported empirically. When every \(n_{i} = 1\), \(a_{j} = m - j + 1\) and \(Z\) is the factorial-base representation of a permutation, so MP-BAT with all multiplicities equal to one enumerates the same set as the permutation BAT of {[}6{]} and advances its counter at the same cost. The two are nevertheless different bijections: the permutation BAT decodes by shifting, and its digit prefixes do not fix route prefixes, so Proposition 3 does not transfer to it.

A lexicographic counter that skips pruned blocks visits the same prefixes, in the same order, as a depth-first branch and bound over the prefix tree with children taken in IC order. The counter formulation is what supplies the bijection of Proposition 2, the completeness argument and the counter windows used for parallelisation; it is not a different search order, and Section 6.3 measures it against a different one. Nor is MP-BAT claimed to be the fastest exact method for the TR-UDSP in general. A Held--Karp dynamic programme {[}28{]} over states (remaining-count vector, last IC) has \(m\prod_{i}\left( n_{i} + 1 \right)\) states, far fewer than \(|\Omega|\) when multiplicities are large; but here the future of a prefix depends on the last completion time of every started IC, through the flow reliability of the next gap in both variants and through the waits in the hard-dwell variant, so the state must carry those times, and with three objectives each state holds a Pareto set rather than a scalar. MP-BAT trades a larger number of expanded prefixes for a traversal memory of \(O(|N|)\) states, to which only the archive of \(|\Phi|\) routes is added, and for an enumeration whose correctness does not depend on how the objectives are defined: schedule-dependent reliability, dwell rules and further objectives change \(Extend\) and \(Bound\) in Algorithm 1 and nothing else.

\hypertarget{computational-results}{%
\section{6. Computational Results}\label{computational-results}}

\hypertarget{instances-parameters-and-implementation}{%
\subsection{6.1. Instances, parameters and implementation}\label{instances-parameters-and-implementation}}

Three synthetic instances and three instances derived from the BeagleBone Black Rev C layout are used (Table 6). Instance 1 is the four-IC board of Figure 1(a): A (20 × 20 mm, line (30, 60)--(30, 80), \(\tau = 15\) s, amounts 10, 10, 10, 8, 7 mg), B (15 × 15, (90, 62.5)--(90, 77.5), \(\tau = 12\), amounts 10, 9, 8, 7), C (10 × 10, (30, 15)--(30, 25), \(\tau = 8\), amounts 5, 4, 3) and D (10 × 10, (90, 15)--(90, 25), \(\tau = 8\), amounts 5, 5, 3). Instance 2 adds E (12 × 12 at (55, 38)--(67, 50), \(\tau = 10\), amounts 6, 5, 4) and gives every IC three passes; Instance 3 gives A and B four passes. The BBB instances take the six top-side BGA/QFN devices of the BeagleBone Black assembly XY report {[}29{]} (Figure 1(b)); pass counts, amounts and flow times are assigned by die area and are part of the instance files; BBB-14, BBB-16 and BBB-20 share the six dies, their positions and flow times, and differ in the pass counts of Table 6 and hence in the pass amounts.

Parameters throughout: \(v = 50\) mm/s, \(\rho = 2\) mg/s, \(\sigma(t) = 0.10\, t\), \(\sigma(p) = 0.05\, p\), \(r^{P}(q) = 1 - 0.002\, q\), \(r_{0} = 0.99\), \(\lambda = 0.001\ s^{- 1}\). MP-BAT is implemented in C in two forms: a serial program and an OpenMP program that distributes the counter windows of Section 5.2 over threads (window depth 6 for up to 17 passes, 8 for BBB-20). All reference fronts were produced on an NVIDIA DGX Spark (20-core Grace CPU, 10 Cortex-X925 and 10 Cortex-A725, aarch64, 121 GB RAM, gcc 13.3, -O3, 20 threads); the GPU is not used.

Three heuristics share the instance reader and evaluation code with MP-BAT (NSGA-II and MSSO test dominance by exact comparison of the stored floating-point objectives, whereas MP-BAT and the seeding routine use the tolerances of Section 5.2): random-restart neighbourhood sampling (PLS: random feasible routes and their order-preserving one-swap neighbours are evaluated and the non-dominated candidates retained), used only to produce warm-start archives; NSGA-II {[}17{]} (population 100, binary tournament on rank and crowding distance, order crossover adapted to multiset strings with probability 0.9, swap mutation with probability 0.3); and a multi-objective simplified swarm optimisation (MSSO) in the form of {[}18,19{]} (population 100, per-particle update from the current position, the personal best or a random member of the external archive with \(C_{g},C_{p},C_{w} = 0.4,0.7,0.9\), the standard setting of {[}21{]}, otherwise a random route). Both population methods keep an external archive of every non-dominated route found, and the multiset order crossover guarantees feasibility, so no repair is needed. Neither was tuned. The supplement documents the historical PLS implementation, which was written as an iterative Pareto local search but, owing to an ordering defect found in review, never accepted a move, and provides a corrected iterative variant that was not used for any reported result; a warm start requires only feasible seed routes, which the sampling supplies.

Six controls are applied, following the practice that an exact method should carry its own verification: (i) without pruning the number of leaves evaluated must equal \(|\Omega|\) of Equation (1), which held for all five instances on which full enumeration was run; (ii) the fronts of the serial and parallel programs must coincide as sets of objective points \(\left( T,\sigma^{*},R \right)\) rounded to four decimals, which held for Instances 1, 2 and BBB-14, a regression check (Section 5.2); (iii) a front is reported as a reference front only if its run completed without hitting a time limit; a reference front is the output of a completed double-precision search with the dominance tolerances of Section 5.2, checked against full enumeration where that was feasible and against the controls below, not an interval-arithmetic certificate; and (iv) an independent single-objective control: a CP-SAT {[}7{]} model of the minimum-\(T\) hard-dwell route (circuit constraint over the pass nodes, time propagation on chosen arcs, and the dwell lags as minimum time lags) was solved to proven optimality on Instances 1, 2 and BBB-14, giving \(T = 87.551\), \(65.818\) and \(56.078\) s, equal to the minimum-\(T\) members of the reference fronts within 1 ms (the CP-SAT model rounds every duration to 1 ms, so its optimum can differ from the exact one by at most 0.5 ms per arc; the solver returned tied routes with the same objective on two of the three); (v) a seed-independence control for warm-started runs: seed routes are re-evaluated by the parallel solver and their stated objectives are ignored, and a run of Instance 1 seeded with an archive whose objectives had been deliberately rounded to two decimals reproduced the reference front with no point missing and none invented (the serial program trusted the stated objectives until v0.4 and returns 121 points on the same rounded seed; every warm-started result reported here was produced by the parallel program, and the single-core development runs of Section 6.3 were seeded with objectives computed by the same evaluator); and (vi) a regression instance for the bound. An earlier version of the implementation computed \(t_{\min}\) over moves between ICs only, omitting the trips from and to \(h\). On the six benchmarks those trips are never the shortest move (0.67 and 0.78 s against 0.20 s on Instances 1--3; 0.47 and 0.54 s against 0.053 s on the BBB boards), so the bound was valid and every reported front is unaffected; the corrected and the earlier definition give identical prefix, leaf and pruning counts on Instances 1, 2 and BBB-14. On a constructed three-IC instance in which one dispense line ends next to \(h\), the earlier bound pruned the minimum-\(T\) route (12.671 s) and returned a front of one point where full enumeration returns two; the corrected bound returns both. The instance is supplied as a regression test, and the solvers of the supplementary material carry the correction (mpbat\_omp v0.5, mpbat\_serial v0.3, mpbat\_lab v0.2).

Because two different routes can have identical objectives (they differ by swapping equal-amount passes whose gaps and moves mirror each other), fronts are reported as sets of objective points; the serial and parallel programs then keep different representatives of such ties, which is why (ii) compares points rather than routes. Fronts are counted as objective points rounded to four decimals; in the reference front files only one pair of rows coincides at that precision (BBB-20, hard-dwell: equal \(T\) and \(\sigma^{*}\), \(R\) differing by \(2.4 \times 10^{-5}\)), so that front of 2,045 routes counts as 2,044 points. The experiments were run as three campaigns on the same machine, each recorded with its repository commit in the campaign log: the exact fronts and their controls, the heuristic runs, and the re-computation of the warm-started fronts with control (v) together with the scoring of the heuristic archives. The validation experiments added in this revision (matched baselines in Section 6.3, Sections 6.5 and 6.6) use an instrumented serial build of the same evaluator, bound and dominance test (mpbat\_lab) and Python scripts, run on one core of an x86-64 cloud container (gcc 13.3, -O3); that build reproduced the reference Instance 1 fronts point for point.

\textbf{Table 6.} Instances. \(|\Omega|\) is the number of feasible routes.

\begin{longtable}[]{@{}
  >{\raggedright\arraybackslash}p{(\columnwidth - 10\tabcolsep) * \real{0.1805}}
  >{\raggedright\arraybackslash}p{(\columnwidth - 10\tabcolsep) * \real{0.0565}}
  >{\raggedright\arraybackslash}p{(\columnwidth - 10\tabcolsep) * \real{0.1280}}
  >{\raggedright\arraybackslash}p{(\columnwidth - 10\tabcolsep) * \real{0.0745}}
  >{\raggedright\arraybackslash}p{(\columnwidth - 10\tabcolsep) * \real{0.1303}}
  >{\raggedright\arraybackslash}p{(\columnwidth - 10\tabcolsep) * \real{0.4302}}@{}}
\toprule\noalign{}
\begin{minipage}[b]{\linewidth}\raggedright
Instance
\end{minipage} & \begin{minipage}[b]{\linewidth}\raggedright
ICs
\end{minipage} & \begin{minipage}[b]{\linewidth}\raggedright
Passes \(n_{i}\)
\end{minipage} & \begin{minipage}[b]{\linewidth}\raggedright
\(|N|\)
\end{minipage} & \begin{minipage}[b]{\linewidth}\raggedright
\(|\Omega|\)
\end{minipage} & \begin{minipage}[b]{\linewidth}\raggedright
Status in this paper
\end{minipage} \\
\midrule\noalign{}
\endhead
\bottomrule\noalign{}
\endlastfoot
Inst-1 & 4 & 5, 4, 3, 3 & 15 & \(1.26 \times 10^{7}\) & reference front, both variants \\
Inst-2 & 5 & 3, 3, 3, 3, 3 & 15 & \(1.68 \times 10^{8}\) & reference front, both variants \\
BBB-14 & 6 & 3, 3, 2, 2, 2, 2 & 14 & \(1.51 \times 10^{8}\) & reference front, both variants \\
Inst-3 & 5 & 4, 4, 3, 3, 3 & 17 & \(2.86 \times 10^{9}\) & reference front, both variants \\
BBB-16 & 6 & 3, 3, 3, 3, 2, 2 & 16 & \(4.04 \times 10^{9}\) & reference front, both variants \\
BBB-20 & 6 & 5, 4, 4, 3, 2, 2 & 20 & \(1.47 \times 10^{12}\) & reference front, both variants (warm-started) \\
\end{longtable}

\hypertarget{exact-pareto-fronts-and-cost-scalarisation-on-instance-1}{%
\subsection{6.2. Exact Pareto fronts and cost scalarisation on Instance 1}\label{exact-pareto-fronts-and-cost-scalarisation-on-instance-1}}

MP-BAT enumerates the 12,612,600 routes of Instance 1 in 0.10 s on 20 cores (1.1 s serial); the hard-dwell front has 194 members and the soft-dwell front 115. Table 7 lists the extreme members and the IC-by-IC baseline; Figure 2 shows the corresponding schedules and Figure 3 the \(T\)--\(R\) projection of both fronts.

\textbf{Table 7.} Extreme members of the Instance 1 fronts and the baseline route. Four hard-dwell routes share the minimum \(T\) = 87.55 s, all Pareto-optimal; the two listed are the extremes, with the smallest \(\sigma^{*}\) and the largest \(R\). \(K_{P} = 0.811\).

\begin{longtable}[]{@{}
  >{\raggedright\arraybackslash}p{(\columnwidth - 12\tabcolsep) * \real{0.1731}}
  >{\raggedright\arraybackslash}p{(\columnwidth - 12\tabcolsep) * \real{0.2166}}
  >{\raggedright\arraybackslash}p{(\columnwidth - 12\tabcolsep) * \real{0.2429}}
  >{\raggedright\arraybackslash}p{(\columnwidth - 12\tabcolsep) * \real{0.0932}}
  >{\raggedright\arraybackslash}p{(\columnwidth - 12\tabcolsep) * \real{0.1117}}
  >{\raggedright\arraybackslash}p{(\columnwidth - 12\tabcolsep) * \real{0.0670}}
  >{\raggedright\arraybackslash}p{(\columnwidth - 12\tabcolsep) * \real{0.0955}}@{}}
\toprule\noalign{}
\begin{minipage}[b]{\linewidth}\raggedright
Variant
\end{minipage} & \begin{minipage}[b]{\linewidth}\raggedright
Criterion
\end{minipage} & \begin{minipage}[b]{\linewidth}\raggedright
Route
\end{minipage} & \begin{minipage}[b]{\linewidth}\raggedright
\(T\) (s)
\end{minipage} & \begin{minipage}[b]{\linewidth}\raggedright
\(\sqrt{\sigma^{*}}\) (s)
\end{minipage} & \begin{minipage}[b]{\linewidth}\raggedright
\(R\)
\end{minipage} & \begin{minipage}[b]{\linewidth}\raggedright
\(R/K_{P}\)
\end{minipage} \\
\midrule\noalign{}
\endhead
\bottomrule\noalign{}
\endlastfoot
hard & min \(T\) & ABCADBCADBACDBA & 87.6 & 0.934 & 0.760 & 0.937 \\
hard & min \(T\) (tie) & ABABDCABDCABDCA & 87.6 & 0.910 & 0.758 & 0.934 \\
hard & min \(\sigma^{*}\) & CAAAAABBBBDDDCC & 182.3 & 0.818 & 0.727 & 0.896 \\
hard & max \(R\) & ADCBDCAABBAABDC & 114.3 & 0.896 & 0.771 & 0.950 \\
soft & min \(T\) & CCCDDDBBBBAAAAA & 65.3 & 0.823 & ≈ 0 & ≈ 0 \\
soft & min \(\sigma^{*}\) & CCAAAAABBBBDDDC & 65.4 & 0.818 & ≈ 0 & ≈ 0 \\
soft & max \(R\) & ABDACBADCBADCBA & 78.6 & 0.954 & 0.720 & 0.888 \\
hard & IC-by-IC baseline & AAAAABBBBCCCDDD & 191.8 & 0.849 & 0.718 & 0.886 \\
\end{longtable}

\includegraphics[width=6in,height=4in]{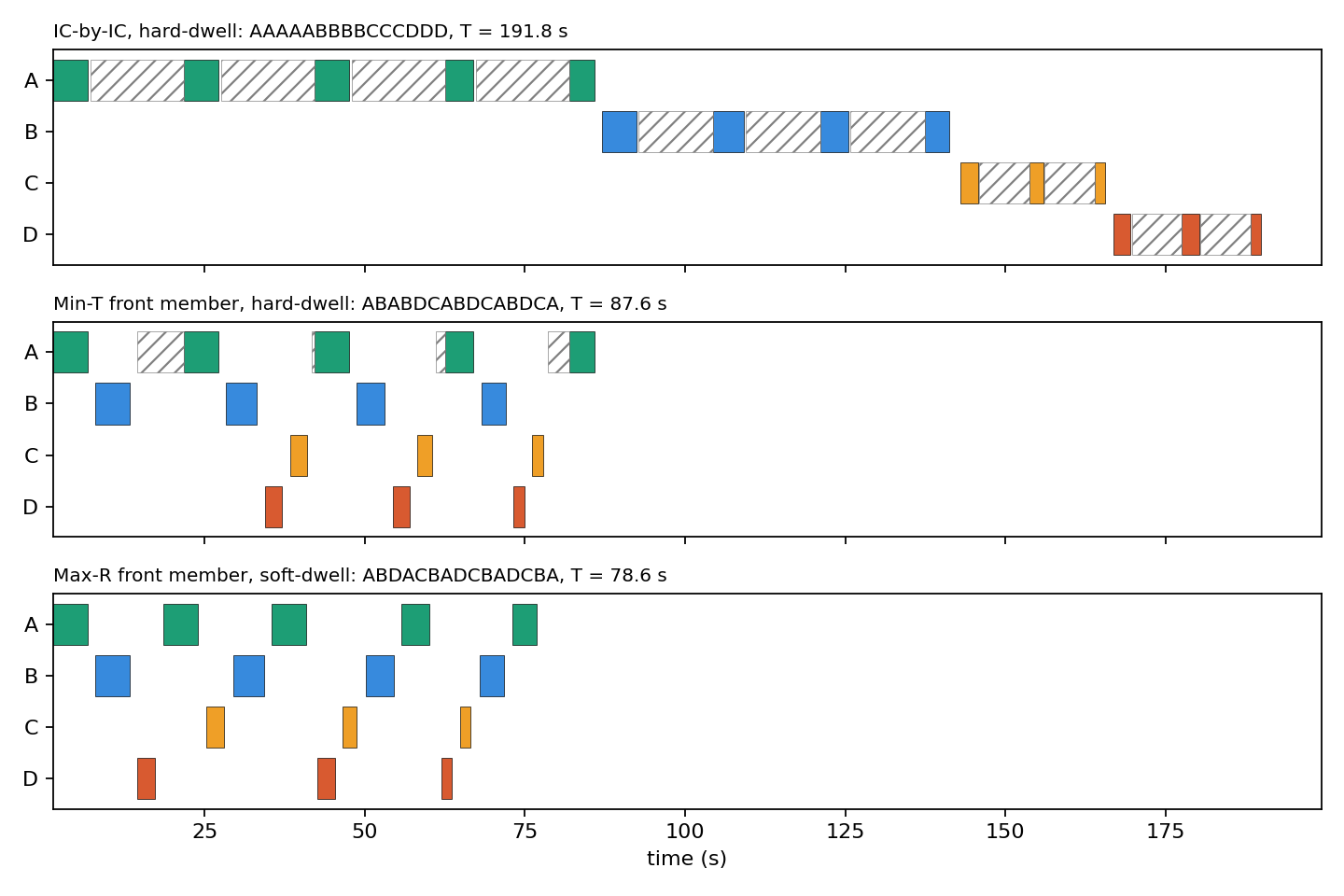}

Figure 2. Schedules of three Instance 1 routes: the IC-by-IC baseline (hard-dwell), the min-\(T\) member of the hard-dwell front and the max-\(R\) member of the soft-dwell front. Bars are passes; hatched bars are forced waiting.

\includegraphics[width=4.5in,height=3in]{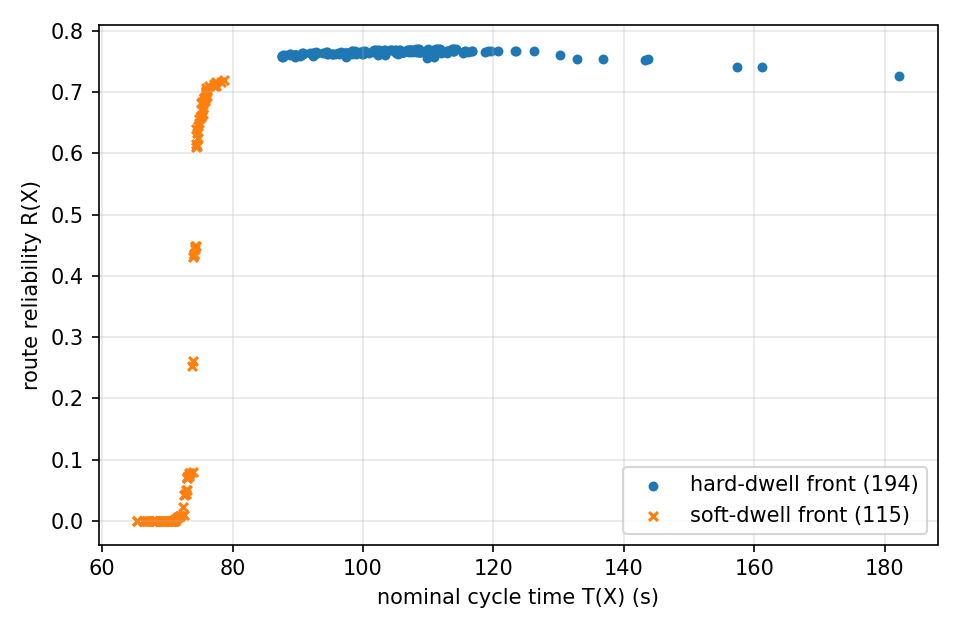}

Figure 3. Nominal cycle time \(T(X)\) against route reliability \(R(X)\) for the hard-dwell (circles) and soft-dwell (crosses) Pareto fronts of Instance 1.

Four observations follow. First, against IC-by-IC processing the hard-dwell min-\(T\) routes cut cycle time by 54 \% and raise reliability. Of the 128 s of nominal dwell, the baseline spends 124.7 s waiting; the min-\(T\) routes wait 11.2--12.4 s and hide the rest behind passes on other ICs. Second, the soft-dwell min-\(T\) and min-\(\sigma^{*}\) members have \(R \approx 0\): they are Pareto-optimal only because nothing beats them on time. A reliability requirement excludes them as surely as the cost model does (Table 9); what an objective adds over a constraint is that the whole trade-off is available before a threshold is chosen. Third, the most reliable soft-dwell route is 8.9 s faster than the most reliable min-\(T\) hard-dwell route at a cost of 0.040 in \(R\), the trade that the cost model prices. Fourth, the pass factor \(K_{P} = 0.811\) is common to every route, so no sequence can exceed \(R = 0.811\) on this instance, and a requirement such as \(R \geq 0.9\) is infeasible whatever the algorithm; the sequencing effect is \(R/K_{P}\), which ranges from 0.886 at the IC-by-IC baseline to 0.950 at the max-\(R\) route.

With \(c_{t} = 0.02\) \$/s (a 72 \$/h line) and \(c_{w} = 0\) (waiting is charged through \(T\)), Table 8 gives the route selected from each front as the failure cost \(c_{f}\) varies. The soft-dwell front supplies the cheaper route only when a failed board costs less than about 231 s of line time (\(c_{f} < 4.63\) \$ at \(c_{t} = 0.02\) \$/s, the exact crossover of the two cost curves); above that the selection moves steadily from the fast end of the hard-dwell front (87.6 s) to its reliable end (114.3 s). Against the baseline at \(c_{f} = 30\) \$ (\(\mathcal{C} = 12.29\) \$/board), the selected route saves 3.39 \$/board and raises throughput from 18.8 to 39.6 boards per hour.

\textbf{Table 8.} Route selected by Equation (9) on Instance 1 as a function of the failure cost \(c_{f}\) (\(c_{t} = 0.02\) \$/s): (a) hard-dwell front, (b) soft-dwell front.

\begin{enumerate}
\def\labelenumi{(\alph{enumi})}
\tightlist
\item
  Hard-dwell front
\end{enumerate}

\begin{longtable}[]{@{}
  >{\raggedright\arraybackslash}p{(\columnwidth - 8\tabcolsep) * \real{0.2018}}
  >{\raggedright\arraybackslash}p{(\columnwidth - 8\tabcolsep) * \real{0.3511}}
  >{\raggedright\arraybackslash}p{(\columnwidth - 8\tabcolsep) * \real{0.1348}}
  >{\raggedright\arraybackslash}p{(\columnwidth - 8\tabcolsep) * \real{0.0969}}
  >{\raggedright\arraybackslash}p{(\columnwidth - 8\tabcolsep) * \real{0.2155}}@{}}
\toprule\noalign{}
\begin{minipage}[b]{\linewidth}\raggedright
\(c_{f}\) (\$)
\end{minipage} & \begin{minipage}[b]{\linewidth}\raggedright
Selected route
\end{minipage} & \begin{minipage}[b]{\linewidth}\raggedright
\(T\) (s)
\end{minipage} & \begin{minipage}[b]{\linewidth}\raggedright
\(R\)
\end{minipage} & \begin{minipage}[b]{\linewidth}\raggedright
\(\mathcal{C}\) (\$/board)
\end{minipage} \\
\midrule\noalign{}
\endhead
\bottomrule\noalign{}
\endlastfoot
2 & ABCADBCADBACDBA & 87.6 & 0.760 & 2.23 \\
5 & ABCADBCADBACDBA & 87.6 & 0.760 & 2.95 \\
10 & ABDCABDCABACDBA & 87.7 & 0.760 & 4.15 \\
30 & ABDACBDCABACDBA & 90.8 & 0.764 & 8.90 \\
100 & ACDBCADCBDAABBA & 101.7 & 0.769 & 25.14 \\
200 & ADCBDCAABBAABDC & 114.3 & 0.771 & 48.16 \\
\end{longtable}

\begin{enumerate}
\def\labelenumi{(\alph{enumi})}
\setcounter{enumi}{1}
\tightlist
\item
  Soft-dwell front
\end{enumerate}

\begin{longtable}[]{@{}
  >{\raggedright\arraybackslash}p{(\columnwidth - 8\tabcolsep) * \real{0.2018}}
  >{\raggedright\arraybackslash}p{(\columnwidth - 8\tabcolsep) * \real{0.3511}}
  >{\raggedright\arraybackslash}p{(\columnwidth - 8\tabcolsep) * \real{0.1348}}
  >{\raggedright\arraybackslash}p{(\columnwidth - 8\tabcolsep) * \real{0.0969}}
  >{\raggedright\arraybackslash}p{(\columnwidth - 8\tabcolsep) * \real{0.2155}}@{}}
\toprule\noalign{}
\begin{minipage}[b]{\linewidth}\raggedright
\(c_{f}\) (\$)
\end{minipage} & \begin{minipage}[b]{\linewidth}\raggedright
Selected route
\end{minipage} & \begin{minipage}[b]{\linewidth}\raggedright
\(T\) (s)
\end{minipage} & \begin{minipage}[b]{\linewidth}\raggedright
\(R\)
\end{minipage} & \begin{minipage}[b]{\linewidth}\raggedright
\(\mathcal{C}\) (\$/board)
\end{minipage} \\
\midrule\noalign{}
\endhead
\bottomrule\noalign{}
\endlastfoot
2 & ABCADBACDBACDBA & 75.9 & 0.707 & 2.11 \\
5 & ABDACBACDBADCBA & 77.5 & 0.716 & 2.97 \\
10 & ABDACBADCBADCBA & 78.6 & 0.720 & 4.37 \\
30 & ABDACBADCBADCBA & 78.6 & 0.720 & 9.97 \\
100 & ABDACBADCBADCBA & 78.6 & 0.720 & 29.57 \\
200 & ABDACBADCBADCBA & 78.6 & 0.720 & 57.57 \\
\end{longtable}

Where a failure cost is hard to state, a reliability requirement selects instead. Table 9 gives the fastest route on the union of both Instance 1 fronts for a range of requirements \(R_{\min}\), read directly off the reference fronts. Up to \(R_{\min} = 0.72\) the soft-dwell front supplies the route (75.9--78.6 s); from 0.74 to 0.76 a hard-dwell route at 87.6--87.7 s; the last percent of attainable reliability costs another 26 s (0.770 at 113.6 s), and nothing above 0.771 is attainable.

\textbf{Table 9.} Requirement-based selection on Instance 1: fastest route with \(R \geq R_{\min}\) over both fronts.

\begin{longtable}[]{@{}
  >{\raggedright\arraybackslash}p{(\columnwidth - 12\tabcolsep) * \real{0.1922}}
  >{\raggedright\arraybackslash}p{(\columnwidth - 12\tabcolsep) * \real{0.1224}}
  >{\raggedright\arraybackslash}p{(\columnwidth - 12\tabcolsep) * \real{0.1134}}
  >{\raggedright\arraybackslash}p{(\columnwidth - 12\tabcolsep) * \real{0.2787}}
  >{\raggedright\arraybackslash}p{(\columnwidth - 12\tabcolsep) * \real{0.1070}}
  >{\raggedright\arraybackslash}p{(\columnwidth - 12\tabcolsep) * \real{0.0769}}
  >{\raggedright\arraybackslash}p{(\columnwidth - 12\tabcolsep) * \real{0.1095}}@{}}
\toprule\noalign{}
\begin{minipage}[b]{\linewidth}\raggedright
\(R_{\min}\)
\end{minipage} & \begin{minipage}[b]{\linewidth}\raggedright
\(R_{\min}/K_{P}\)
\end{minipage} & \begin{minipage}[b]{\linewidth}\raggedright
Variant
\end{minipage} & \begin{minipage}[b]{\linewidth}\raggedright
Route
\end{minipage} & \begin{minipage}[b]{\linewidth}\raggedright
\(T\) (s)
\end{minipage} & \begin{minipage}[b]{\linewidth}\raggedright
\(R\)
\end{minipage} & \begin{minipage}[b]{\linewidth}\raggedright
\(R/K_{P}\)
\end{minipage} \\
\midrule\noalign{}
\endhead
\bottomrule\noalign{}
\endlastfoot
0.70 & 0.863 & soft & ABCADBACDBACDBA & 75.9 & 0.707 & 0.872 \\
0.72 & 0.888 & soft & ABDACBADCBADCBA & 78.6 & 0.720 & 0.888 \\
0.74 & 0.913 & hard & ABABDCABDCABDCA & 87.6 & 0.758 & 0.934 \\
0.76 & 0.937 & hard & ABDCABDCABACDBA & 87.7 & 0.760 & 0.937 \\
0.765 & 0.943 & hard & ABADBCDABACDBAC & 94.4 & 0.766 & 0.945 \\
0.77 & 0.950 & hard & CADBCDAABBAABDC & 113.6 & 0.770 & 0.950 \\
\end{longtable}

\hypertarget{enumeration-efficiency-and-larger-boards}{%
\subsection{6.3. Enumeration efficiency and larger boards}\label{enumeration-efficiency-and-larger-boards}}

Table 10 reports, for every instance, the number of vectors the alternatives of Section 5.1 would generate, the leaves MP-BAT evaluates with and without the prefix bound, and the wall time on 20 cores. Full enumeration runs at \(1.2\)--\(1.4 \times 10^{8}\) leaves per second and its leaf count equals \(|\Omega|\) exactly on every instance, which is the completeness check (i) of Section 6.1. The prefix bound removes 98.2 \% of the leaves on Instance 3 and 99.5 \% on BBB-16 in the hard-dwell variant and about 90 \% in the soft-dwell variant, whose bound is weaker because no dwell is forced; the wall-time gain is smaller than the leaf reduction because pruned prefixes still cost a dominance test each. Figure 4 shows wall time against \(|\Omega|\) for both modes.

\textbf{Table 10.} Vectors generated, leaves evaluated and wall time (20 cores; serial times in parentheses). Pruned leaf counts are means over five repeats; they vary by about 1 \% between runs because the parallel archives grow in a thread-dependent order.

\begin{longtable}[]{@{}
  >{\raggedright\arraybackslash}p{(\columnwidth - 10\tabcolsep) * \real{0.1591}}
  >{\raggedright\arraybackslash}p{(\columnwidth - 10\tabcolsep) * \real{0.0962}}
  >{\raggedright\arraybackslash}p{(\columnwidth - 10\tabcolsep) * \real{0.3053}}
  >{\raggedright\arraybackslash}p{(\columnwidth - 10\tabcolsep) * \real{0.2291}}
  >{\raggedright\arraybackslash}p{(\columnwidth - 10\tabcolsep) * \real{0.0745}}
  >{\raggedright\arraybackslash}p{(\columnwidth - 10\tabcolsep) * \real{0.1359}}@{}}
\toprule\noalign{}
\begin{minipage}[b]{\linewidth}\raggedright
Instance
\end{minipage} & \begin{minipage}[b]{\linewidth}\raggedright
Variant
\end{minipage} & \begin{minipage}[b]{\linewidth}\raggedright
Method
\end{minipage} & \begin{minipage}[b]{\linewidth}\raggedright
Vectors / leaves evaluated
\end{minipage} & \begin{minipage}[b]{\linewidth}\raggedright
Front
\end{minipage} & \begin{minipage}[b]{\linewidth}\raggedright
Time (s)
\end{minipage} \\
\midrule\noalign{}
\endhead
\bottomrule\noalign{}
\endlastfoot
Inst-1 & hard & fixed-radix BAT, filtered & \(1.07 \times 10^{9}\) / 12,612,600 & 194 & --- \\
Inst-1 & hard & permutation BAT on pass nodes, filtered & \(1.3 \times 10^{12}\) / 12,612,600 & 194 & --- \\
Inst-1 & hard & MP-BAT, full & 12,612,600 & 194 & 0.10 \\
Inst-1 & hard & MP-BAT + bound & 1,832,730 (1,533,880 serial) & 194 & 0.11 (1.13) \\
Inst-2 & hard & MP-BAT, full & 168,168,000 & 158 & 1.16 \\
Inst-2 & hard & MP-BAT + bound & 6,758,435 (4,370,667 serial) & 158 & 0.56 (4.00) \\
BBB-14 & hard & MP-BAT, full & 151,351,200 & 255 & 1.10 \\
BBB-14 & hard & MP-BAT + bound & 4,378,392 (1,875,554 serial) & 255 & 0.43 (3.41) \\
Inst-3 & hard & MP-BAT, full & 2,858,856,000 & 349 & 21.01~±~0.64 \\
Inst-3 & hard & MP-BAT + bound & 50,536,761 & 349 & 9.45~±~0.17 \\
Inst-3 & soft & MP-BAT + bound & 318,933,867 & 287 & 30.69~±~0.36 \\
BBB-16 & hard & MP-BAT, full & 4,036,032,000 & 342 & 29.00~±~0.18 \\
BBB-16 & hard & MP-BAT + bound & 22,051,493 & 342 & 5.05~±~0.14 \\
BBB-16 & soft & MP-BAT + bound & 48,974,504 & 191 & 6.44~±~0.24 \\
\end{longtable}

\includegraphics[width=4.8in,height=3.36in]{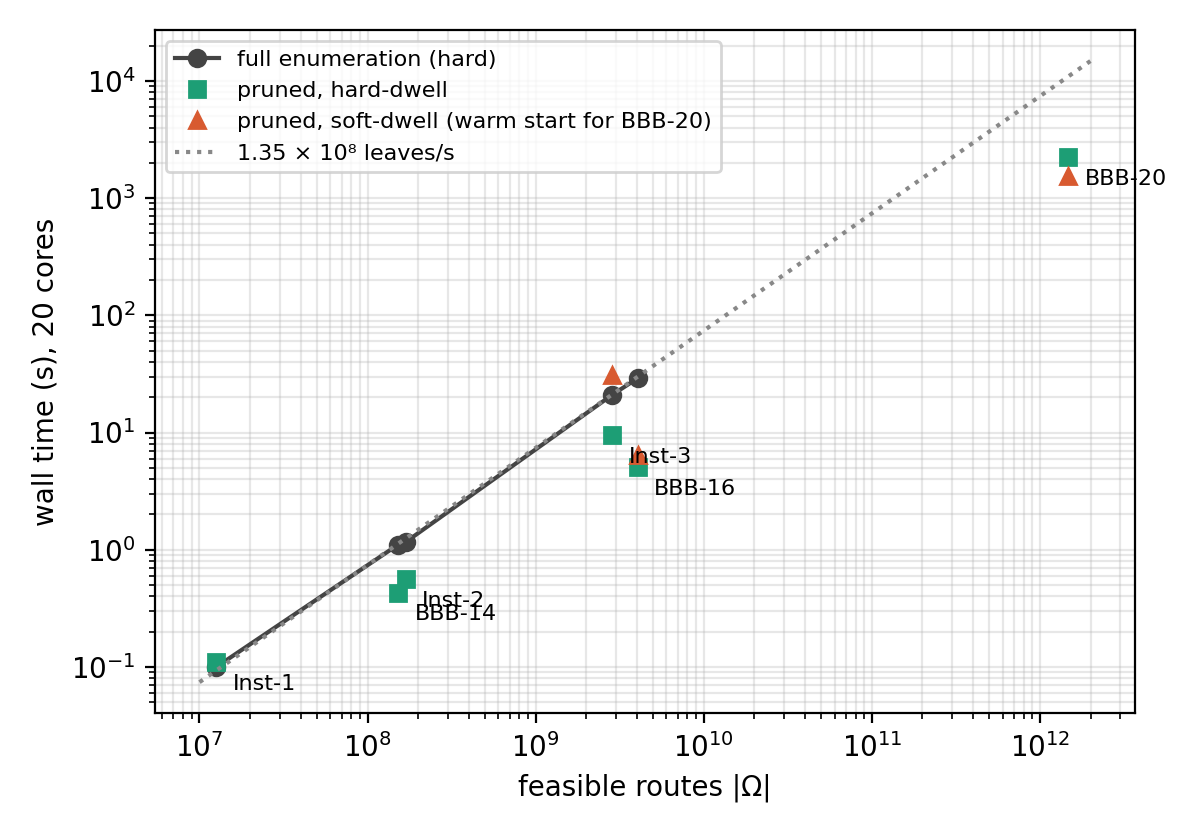}

Figure 4. Wall time on 20 cores against the number of feasible routes \(|\Omega|\): full enumeration (line), pruned hard-dwell and pruned soft-dwell runs (markers), and the \(1.35 \times 10^{8}\) leaves/s reference. BBB-20 was warm-started from the PLS front.

On BBB-14 and BBB-16 the max-\(R\) route is the same under both variants and the soft front adds only the \(R \approx 0\) fast routes: on these boards six ICs provide enough other work to hide every flow time, so the dwell conflict is weaker than on Instance 1; the number of ICs alone does not guarantee this, which depends on the flow times relative to the intervening work.

BBB-20 has \(1.47 \times 10^{12}\) feasible routes, about 3 h of full enumeration on this machine at the measured rate and more than the 2 h allotted per run. Two observations from single-core development runs motivate the warm start: a cold-started pruned MP-BAT given 240 s spent its whole budget in the A-first corner of the tree and its 844-member archive was dominated in full by a 0.4 s PLS front, whereas the same budget seeded with that front pruned \(1.1 \times 10^{9}\) prefixes and found 254 new members. On 20 cores with a 2 h budget per run, MP-BAT completed both BBB-20 searches from a PLS seed and obtained the reference fronts (Table 11): the hard-dwell front, 2,044 objective points, in 38 min after pruning \(5.4 \times 10^{10}\) prefixes and evaluating \(5.6 \times 10^{8}\) leaves; the soft-dwell front, 514 points, in 25 min. Not one of the 106 hard-dwell or 29 soft-dwell PLS seed members survived; the exact hard-dwell front dominates the seed at both ends, starting at \(T = 85.9\) s against the seed's best 89.2 s and reaching \(R = 0.7359\) against 0.7351. A third run seeded with the union of thirty NSGA-II archives (7,943 routes, Section 6.4) returned exactly the same 2,044 points, which is the seed-independence the method promises, and it shows what a seed is worth: the stronger seed cut the evaluated leaves by 75 \% but lengthened the run by 8 \%, because every prefix bound test now scans an archive seeded with thousands of points. A seed of a hundred routes from a cheap heuristic is enough to make the bound bite from the first vector; a larger one buys pruning that its own dominance tests consume.

\textbf{Table 11.} BBB-20: reference fronts by warm-started parallel MP-BAT (20 cores, 1,173,018 counter windows at depth 8, seed routes re-evaluated exactly). All three runs reproduce the fronts of the first campaign with no point missing and none new. ``Seed surviving'' counts seed routes that are members of the final front, that is, Pareto-optimal routes the heuristic had already found.

\begin{longtable}[]{@{}
  >{\raggedright\arraybackslash}p{(\columnwidth - 12\tabcolsep) * \real{0.1369}}
  >{\raggedright\arraybackslash}p{(\columnwidth - 12\tabcolsep) * \real{0.1967}}
  >{\raggedright\arraybackslash}p{(\columnwidth - 12\tabcolsep) * \real{0.1561}}
  >{\raggedright\arraybackslash}p{(\columnwidth - 12\tabcolsep) * \real{0.1452}}
  >{\raggedright\arraybackslash}p{(\columnwidth - 12\tabcolsep) * \real{0.1173}}
  >{\raggedright\arraybackslash}p{(\columnwidth - 12\tabcolsep) * \real{0.0850}}
  >{\raggedright\arraybackslash}p{(\columnwidth - 12\tabcolsep) * \real{0.1628}}@{}}
\toprule\noalign{}
\begin{minipage}[b]{\linewidth}\raggedright
Variant
\end{minipage} & \begin{minipage}[b]{\linewidth}\raggedright
Seed (routes)
\end{minipage} & \begin{minipage}[b]{\linewidth}\raggedright
Leaves evaluated
\end{minipage} & \begin{minipage}[b]{\linewidth}\raggedright
Pruned prefixes
\end{minipage} & \begin{minipage}[b]{\linewidth}\raggedright
Front points
\end{minipage} & \begin{minipage}[b]{\linewidth}\raggedright
Time (s)
\end{minipage} & \begin{minipage}[b]{\linewidth}\raggedright
Seed surviving
\end{minipage} \\
\midrule\noalign{}
\endhead
\bottomrule\noalign{}
\endlastfoot
hard & PLS (106) & \(5.57 \times 10^{8}\) & \(5.44 \times 10^{10}\) & 2044 & 2260 & 0~of~106 \\
hard & NSGA-II~(7,943) & \(1.40 \times 10^{8}\) & \(3.03 \times 10^{10}\) & 2044 & 2440 & 1,541~of~7,943 \\
soft & PLS (29) & \(2.16 \times 10^{9}\) & \(1.25 \times 10^{11}\) & 514 & 1523 & 0~of~29 \\
\end{longtable}

\emph{Matched exact baselines.} Table 12 runs one serial evaluator in five configurations on Instances 1, 2 and BBB-14: full enumeration; the counter (depth-first) with the bound, cold and warm-started from one NSGA-II archive of \(10^{6}\) evaluations; and best-first branch and bound on \(T_{lb}\) with the same bound, archive and tolerances, cold and warm. All configurations that completed returned identical fronts. Three results follow. First, the bound reduces evaluated leaves 3--99-fold, but on one core it saves at most a factor of 3.7 in time, because every prefix costs a dominance test (13.9 million tests against 12.6 million leaves on Instance 1 hard); on Instance 1 it saves nothing. Second, the warm start cuts expanded prefixes by 20--47 \% and time by 9--41 \%. Third, best-first expands the same prefixes as the warm-started counter to within 0.1 \%, because its first leaves are the fast routes that make the bound bite, and gains no warm-start benefit of its own; it is 5--15 times slower, needs 252--1,231 MB against the counter's 2 MB, and on the two soft-dwell instances with \(1.5\)--\(1.7 \times 10^{8}\) routes its open set exceeded 12 million prefixes (3 GB) before completion. The counter with a warm start therefore matches best-first in prefixes within 2 MB of memory, which is what makes the \(10^{12}\)-route run possible.

\textbf{Table 12.} Matched exact baselines (one x86-64 core, same evaluator, bound, archive and tolerances). Entries are expanded prefixes (millions) / evaluated leaves (millions) / time (s) for the counter, and expanded prefixes (millions) / peak memory (MB) / time (s) for best-first; warm runs are seeded with one NSGA-II archive of \(10^{6}\) evaluations. Warm best-first runs expanded the same prefixes as cold ones and are omitted.

\begin{longtable}[]{@{}
  >{\raggedright\arraybackslash}p{(\columnwidth - 10\tabcolsep) * \real{0.1555}}
  >{\raggedright\arraybackslash}p{(\columnwidth - 10\tabcolsep) * \real{0.0962}}
  >{\raggedright\arraybackslash}p{(\columnwidth - 10\tabcolsep) * \real{0.1865}}
  >{\raggedright\arraybackslash}p{(\columnwidth - 10\tabcolsep) * \real{0.1689}}
  >{\raggedright\arraybackslash}p{(\columnwidth - 10\tabcolsep) * \real{0.1689}}
  >{\raggedright\arraybackslash}p{(\columnwidth - 10\tabcolsep) * \real{0.2239}}@{}}
\toprule\noalign{}
\begin{minipage}[b]{\linewidth}\raggedright
Instance
\end{minipage} & \begin{minipage}[b]{\linewidth}\raggedright
Variant
\end{minipage} & \begin{minipage}[b]{\linewidth}\raggedright
Full enumeration (s)
\end{minipage} & \begin{minipage}[b]{\linewidth}\raggedright
Counter + bound, cold
\end{minipage} & \begin{minipage}[b]{\linewidth}\raggedright
Counter + bound, warm
\end{minipage} & \begin{minipage}[b]{\linewidth}\raggedright
Best-first + bound, cold
\end{minipage} \\
\midrule\noalign{}
\endhead
\bottomrule\noalign{}
\endlastfoot
Inst-1 & hard & 2.9 & 7.0 / 1.0 / 3.2 & 5.6 / 0.6 / 2.9 & 5.6 / 252 / 15.6 \\
Inst-1 & soft & 3.6 & 15.4 / 4.0 / 3.7 & 11.8 / 2.1 / 2.9 & 11.8 / 1,231 / 42.9 \\
Inst-2 & hard & 27.3 & 30.0 / 3.4 / 12.1 & 19.6 / 1.5 / 8.2 & 19.6 / 1,058 / 89.0 \\
Inst-2 & soft & 58.4 & 106.1 / 19.3 / 43.9 & 74.8 / 9.0 / 32.8 & open set \textgreater{} 12 M (3 GB) \\
BBB-14 & hard & 24.4 & 15.7 / 1.5 / 9.0 & 9.4 / 0.6 / 6.8 & 9.4 / 679 / 43.6 \\
BBB-14 & soft & 39.6 & 36.3 / 4.2 / 10.7 & 19.2 / 1.5 / 6.3 & open set \textgreater{} 12 M (3 GB) \\
\end{longtable}

\hypertarget{metaheuristics-against-the-reference-fronts}{%
\subsection{6.4. Metaheuristics against the reference fronts}\label{metaheuristics-against-the-reference-fronts}}

With a reference front available for every instance--variant pair, the question a practitioner asks can be answered exactly: how much of the front does a metaheuristic find, and how much does it cost to rely on it? NSGA-II and MSSO were run on all twelve instance--variant pairs with budgets of \(10^{6}\) and \(10^{7}\) evaluations and 30 seeds each, 1,440 runs in total (8.6 min on 20 cores; a single \(10^{7}\)-evaluation run takes 5--30 s on one core). Each archive is scored against the reference front by the number of distinct reference points it contains, by the ratio of its hypervolume {[}30{]} to the front's (exact three-dimensional hypervolume, objectives normalised by the front's ranges, reference point 1.1), by the inverted generational distance (IGD) {[}31{]}, and by the gaps of its extreme members to the reference minimum \(T\) and maximum \(R\). Table 13 reports recovery and hypervolume ratio; the full scores, including IGD, are in the supplementary material.

\textbf{Table 13.} Reference points recovered (mean ± SD over 30 seeds) and mean hypervolume ratio at \(10^{6}\) and \(10^{7}\) evaluations: (a) NSGA-II, (b) MSSO.

\begin{enumerate}
\def\labelenumi{(\alph{enumi})}
\tightlist
\item
  NSGA-II
\end{enumerate}

\begin{longtable}[]{@{}
  >{\raggedright\arraybackslash}p{(\columnwidth - 12\tabcolsep) * \real{0.1781}}
  >{\raggedright\arraybackslash}p{(\columnwidth - 12\tabcolsep) * \real{0.0962}}
  >{\raggedright\arraybackslash}p{(\columnwidth - 12\tabcolsep) * \real{0.0745}}
  >{\raggedright\arraybackslash}p{(\columnwidth - 12\tabcolsep) * \real{0.1722}}
  >{\raggedright\arraybackslash}p{(\columnwidth - 12\tabcolsep) * \real{0.1534}}
  >{\raggedright\arraybackslash}p{(\columnwidth - 12\tabcolsep) * \real{0.1722}}
  >{\raggedright\arraybackslash}p{(\columnwidth - 12\tabcolsep) * \real{0.1534}}@{}}
\toprule\noalign{}
\begin{minipage}[b]{\linewidth}\raggedright
Instance
\end{minipage} & \begin{minipage}[b]{\linewidth}\raggedright
Variant
\end{minipage} & \begin{minipage}[b]{\linewidth}\raggedright
Front
\end{minipage} & \begin{minipage}[b]{\linewidth}\raggedright
Recovered, \(10^{6}\)
\end{minipage} & \begin{minipage}[b]{\linewidth}\raggedright
HV ratio, \(10^{6}\)
\end{minipage} & \begin{minipage}[b]{\linewidth}\raggedright
Recovered, \(10^{7}\)
\end{minipage} & \begin{minipage}[b]{\linewidth}\raggedright
HV ratio, \(10^{7}\)
\end{minipage} \\
\midrule\noalign{}
\endhead
\bottomrule\noalign{}
\endlastfoot
Inst-1 & hard & 194 & 176~±~4 & 0.9995 & 193~±~1 & 1.0000 \\
Inst-1 & soft & 115 & 111~±~1 & 0.9985 & 115~±~0 & 0.9999 \\
Inst-2 & hard & 158 & 147~±~4 & 0.9964 & 156~±~2 & 0.9981 \\
Inst-2 & soft & 194 & 185~±~3 & 0.9999 & 193~±~1 & 1.0000 \\
BBB-14 & hard & 255 & 228~±~4 & 0.9993 & 249~±~2 & 0.9998 \\
BBB-14 & soft & 160 & 155~±~7 & 0.9996 & 159~±~0 & 1.0000 \\
Inst-3 & hard & 349 & 268~±~8 & 0.9986 & 323~±~4 & 0.9996 \\
Inst-3 & soft & 287 & 255~±~5 & 0.9997 & 282~±~2 & 1.0000 \\
BBB-16 & hard & 342 & 232~±~12 & 0.9963 & 286~±~15 & 0.9986 \\
BBB-16 & soft & 191 & 168~±~9 & 0.9992 & 177~±~3 & 0.9999 \\
BBB-20 & hard & 2044 & 248~±~62 & 0.9894 & 611~±~107 & 0.9963 \\
BBB-20 & soft & 514 & 229~±~33 & 0.9912 & 342~±~44 & 0.9969 \\
\end{longtable}

\begin{enumerate}
\def\labelenumi{(\alph{enumi})}
\setcounter{enumi}{1}
\tightlist
\item
  MSSO
\end{enumerate}

\begin{longtable}[]{@{}
  >{\raggedright\arraybackslash}p{(\columnwidth - 12\tabcolsep) * \real{0.1781}}
  >{\raggedright\arraybackslash}p{(\columnwidth - 12\tabcolsep) * \real{0.0962}}
  >{\raggedright\arraybackslash}p{(\columnwidth - 12\tabcolsep) * \real{0.0745}}
  >{\raggedright\arraybackslash}p{(\columnwidth - 12\tabcolsep) * \real{0.1722}}
  >{\raggedright\arraybackslash}p{(\columnwidth - 12\tabcolsep) * \real{0.1534}}
  >{\raggedright\arraybackslash}p{(\columnwidth - 12\tabcolsep) * \real{0.1722}}
  >{\raggedright\arraybackslash}p{(\columnwidth - 12\tabcolsep) * \real{0.1534}}@{}}
\toprule\noalign{}
\begin{minipage}[b]{\linewidth}\raggedright
Instance
\end{minipage} & \begin{minipage}[b]{\linewidth}\raggedright
Variant
\end{minipage} & \begin{minipage}[b]{\linewidth}\raggedright
Front
\end{minipage} & \begin{minipage}[b]{\linewidth}\raggedright
Recovered, \(10^{6}\)
\end{minipage} & \begin{minipage}[b]{\linewidth}\raggedright
HV ratio, \(10^{6}\)
\end{minipage} & \begin{minipage}[b]{\linewidth}\raggedright
Recovered, \(10^{7}\)
\end{minipage} & \begin{minipage}[b]{\linewidth}\raggedright
HV ratio, \(10^{7}\)
\end{minipage} \\
\midrule\noalign{}
\endhead
\bottomrule\noalign{}
\endlastfoot
Inst-1 & hard & 194 & 166~±~4 & 0.9994 & 194~±~0 & 1.0000 \\
Inst-1 & soft & 115 & 105~±~3 & 0.9989 & 115~±~0 & 1.0000 \\
Inst-2 & hard & 158 & 135~±~5 & 0.9970 & 156~±~2 & 0.9990 \\
Inst-2 & soft & 194 & 163~±~6 & 0.9991 & 192~±~1 & 1.0000 \\
BBB-14 & hard & 255 & 216~±~7 & 0.9983 & 253~±~1 & 0.9999 \\
BBB-14 & soft & 160 & 147~±~6 & 0.9999 & 160~±~0 & 1.0000 \\
Inst-3 & hard & 349 & 165~±~9 & 0.9933 & 311~±~8 & 0.9984 \\
Inst-3 & soft & 287 & 180~±~8 & 0.9951 & 270~±~3 & 0.9998 \\
BBB-16 & hard & 342 & 187~±~12 & 0.9948 & 317~±~6 & 0.9992 \\
BBB-16 & soft & 191 & 141~±~8 & 0.9975 & 179~±~2 & 1.0000 \\
BBB-20 & hard & 2044 & 66~±~20 & 0.9671 & 647~±~54 & 0.9921 \\
BBB-20 & soft & 514 & 94~±~11 & 0.9702 & 285~±~17 & 0.9962 \\
\end{longtable}

Four things follow from Table 13. First, on the five boards with up to \(4 \times 10^{9}\) routes, \(10^{7}\) evaluations recover 84--100 \% of the reference points: MSSO returns every point of Instance 1 (both variants) and of BBB-14 (soft) in all 30 seeds, and the worst case is 84 \% (NSGA-II on BBB-16, hard). On these boards the exact method is a verification tool; the heuristics find the front. Second, on BBB-20 they do not: at \(10^{7}\) evaluations NSGA-II holds 611 ± 107 of the 2,044 hard-dwell points and MSSO 647 ± 54, about 30 \%, and 67 \% and 55 \% of the 514 soft-dwell points, with a spread across seeds (SD up to 107 points) that a single run would not reveal. Third, the hypervolume ratio does not show this: it is 0.996 for NSGA-II on BBB-20 hard while 70 \% of the front is missing, because the missing points lie in the interior of the front where neighbouring points differ by hundredths of a second and thousandths in \(R\), and the hypervolume is set by the extremes, which every method finds (at \(10^{7}\) evaluations the minimum-\(T\) gap is 0.00 s on all twelve pairs and the maximum-\(R\) gap is at most 0.0006 except for MSSO on BBB-20 soft, 0.0066). A hypervolume ratio above 0.99 is compatible with two thirds of the Pareto-optimal points never having been seen. Most of the missing points differ only in \(\sigma^{*}\): of the 2,044 hard-dwell points of BBB-20, 182 are non-dominated in \((T,R)\) alone, and NSGA-II recovers 26 ± 9 \% and MSSO 11 ± 5 \% of those (soft-dwell: 108 points, 73 ± 9 \% and 57 ± 4 \%); on the five smaller boards both recover 79--100 \% of the \((T,R)\) front. Fourth, which heuristic is better depends on the budget, and Table 14 tests it. Each comparison is a two-sided Mann--Whitney rank-sum test over the 30 independent runs of each method, with the Vargha--Delaney effect size \(\hat{A}_{12}\) (the probability that an MSSO run beats an NSGA-II run; 0.5 is no difference) and Holm correction over the twelve instance--variant pairs. At \(10^{6}\) evaluations NSGA-II recovers more points on all twelve pairs (\(\hat{A}_{12} \leq 0.08\), \(p_{\text{Holm}} < 10^{-7}\)). At \(10^{7}\) evaluations neither dominates: MSSO is significantly better on five pairs (Instance 1 hard and both variants of BBB-14 and BBB-16), NSGA-II on four (Instance 2 soft, both variants of Instance 3, BBB-20 soft) and three show no significant difference; over the twelve pairs the mean recovery fractions do not differ (Wilcoxon signed-rank, \(p = 0.81\)). MSSO's archive-guided update uses the larger budget better on the BeagleBone boards, while NSGA-II's crowding keeps its advantage on the synthetic instances with more passes per die.

\textbf{Table 14.} MSSO against NSGA-II at \(10^{7}\) evaluations: two-sided Mann--Whitney tests over 30 runs of each method, Vargha--Delaney effect size \(\hat{A}_{12}\) (probability that an MSSO run beats an NSGA-II run) and Holm-adjusted \(p\)-values over the twelve pairs, for the number of reference points recovered and for the hypervolume ratio. ``Better'' names the method that is significantly better at \(\alpha = 0.05\); ``---'' means no significant difference.

\begin{longtable}[]{@{}
  >{\raggedright\arraybackslash}p{(\columnwidth - 14\tabcolsep) * \real{0.1565}}
  >{\raggedright\arraybackslash}p{(\columnwidth - 14\tabcolsep) * \real{0.0962}}
  >{\raggedright\arraybackslash}p{(\columnwidth - 14\tabcolsep) * \real{0.1596}}
  >{\raggedright\arraybackslash}p{(\columnwidth - 14\tabcolsep) * \real{0.1064}}
  >{\raggedright\arraybackslash}p{(\columnwidth - 14\tabcolsep) * \real{0.1152}}
  >{\raggedright\arraybackslash}p{(\columnwidth - 14\tabcolsep) * \real{0.1445}}
  >{\raggedright\arraybackslash}p{(\columnwidth - 14\tabcolsep) * \real{0.1064}}
  >{\raggedright\arraybackslash}p{(\columnwidth - 14\tabcolsep) * \real{0.1152}}@{}}
\toprule\noalign{}
\begin{minipage}[b]{\linewidth}\raggedright
Instance
\end{minipage} & \begin{minipage}[b]{\linewidth}\raggedright
Variant
\end{minipage} & \begin{minipage}[b]{\linewidth}\raggedright
\(\hat{A}_{12}\), recovered
\end{minipage} & \begin{minipage}[b]{\linewidth}\raggedright
\(p_{\text{Holm}}\)
\end{minipage} & \begin{minipage}[b]{\linewidth}\raggedright
Better
\end{minipage} & \begin{minipage}[b]{\linewidth}\raggedright
\(\hat{A}_{12}\), HV ratio
\end{minipage} & \begin{minipage}[b]{\linewidth}\raggedright
\(p_{\text{Holm}}\)
\end{minipage} & \begin{minipage}[b]{\linewidth}\raggedright
Better
\end{minipage} \\
\midrule\noalign{}
\endhead
\bottomrule\noalign{}
\endlastfoot
Inst-1 & hard & 0.79 & \textless{} 0.001 & MSSO & 0.80 & \textless{} 0.001 & MSSO \\
Inst-1 & soft & 0.54 & 0.874 & --- & 0.54 & 0.683 & --- \\
Inst-2 & hard & 0.54 & 0.874 & --- & 0.67 & 0.072 & --- \\
Inst-2 & soft & 0.22 & \textless{} 0.001 & NSGA-II & 0.22 & \textless{} 0.001 & NSGA-II \\
BBB-14 & hard & 0.97 & \textless{} 0.001 & MSSO & 0.92 & \textless{} 0.001 & MSSO \\
BBB-14 & soft & 0.82 & \textless{} 0.001 & MSSO & 0.82 & \textless{} 0.001 & MSSO \\
Inst-3 & hard & 0.08 & \textless{} 0.001 & NSGA-II & 0.13 & \textless{} 0.001 & NSGA-II \\
Inst-3 & soft & 0.00 & \textless{} 0.001 & NSGA-II & 0.00 & \textless{} 0.001 & NSGA-II \\
BBB-16 & hard & 0.97 & \textless{} 0.001 & MSSO & 0.75 & 0.004 & MSSO \\
BBB-16 & soft & 0.76 & 0.002 & MSSO & 0.78 & \textless{} 0.001 & MSSO \\
BBB-20 & hard & 0.59 & 0.685 & --- & 0.08 & \textless{} 0.001 & NSGA-II \\
BBB-20 & soft & 0.12 & \textless{} 0.001 & NSGA-II & 0.43 & 0.683 & --- \\
\end{longtable}

A front of two thousand points is not something a line engineer reads; it is something a decision rule is applied to, either by clustering it to a handful of representative routes {[}32{]} or, as here, by scalarising it with the costs the line actually pays. Whether the missing interior matters depends on what the front is used for. Table 15 applies the cost scalarisation of Equation (9) to BBB-20 hard-dwell: for each failure cost it compares the cheapest reference route with the cheapest route in each heuristic archive at \(10^{7}\) evaluations. The regret is small: 0.03--0.23 \% of the cost per board on average, and 0.53 \% in the worst run (NSGA-II, 0.29 \$ at \(c_{f} = 200\) \$), and the reference optimum itself is in the archive in only a minority of seeds (NSGA-II finds the \(c_{f} = 100\) \$ route in 22 of 30 seeds, the \(c_{f} = 2\) \$ route in none). A line engineer using either heuristic would select a route whose excess cost over the best available is below 0.1 \$ per board on average (0.03--0.23 \% of the cost) and 0.29 \$ (0.53 \%) in the worst run, and would not know it. That last clause is the point of the exact method: it does not replace the heuristic on the shop floor, it establishes, once per board design, how far the heuristic's answer can be from the best one.

\textbf{Table 15.} Decision regret on BBB-20 (hard-dwell, \(c_{t} = 0.02\) \$/s): cost per board of the reference selection and the excess cost of the best route in each heuristic archive at \(10^{7}\) evaluations (mean ± SD over 30 seeds). ``Hit'' is the number of seeds whose archive contains the reference selection.

\begin{longtable}[]{@{}
  >{\raggedright\arraybackslash}p{(\columnwidth - 12\tabcolsep) * \real{0.1258}}
  >{\raggedright\arraybackslash}p{(\columnwidth - 12\tabcolsep) * \real{0.1778}}
  >{\raggedright\arraybackslash}p{(\columnwidth - 12\tabcolsep) * \real{0.2174}}
  >{\raggedright\arraybackslash}p{(\columnwidth - 12\tabcolsep) * \real{0.2042}}
  >{\raggedright\arraybackslash}p{(\columnwidth - 12\tabcolsep) * \real{0.0527}}
  >{\raggedright\arraybackslash}p{(\columnwidth - 12\tabcolsep) * \real{0.1694}}
  >{\raggedright\arraybackslash}p{(\columnwidth - 12\tabcolsep) * \real{0.0527}}@{}}
\toprule\noalign{}
\begin{minipage}[b]{\linewidth}\raggedright
\(c_{f}\) (\$)
\end{minipage} & \begin{minipage}[b]{\linewidth}\raggedright
Reference \(\mathcal{C}\) (\$)
\end{minipage} & \begin{minipage}[b]{\linewidth}\raggedright
Selected \(T\) (s), \(R\)
\end{minipage} & \begin{minipage}[b]{\linewidth}\raggedright
NSGA-II regret (\$)
\end{minipage} & \begin{minipage}[b]{\linewidth}\raggedright
Hit
\end{minipage} & \begin{minipage}[b]{\linewidth}\raggedright
MSSO regret (\$)
\end{minipage} & \begin{minipage}[b]{\linewidth}\raggedright
Hit
\end{minipage} \\
\midrule\noalign{}
\endhead
\bottomrule\noalign{}
\endlastfoot
2 & 2.284 & 85.9, 0.7172 & 0.003~±~0.001 & 0 & 0.003~±~0.001 & 0 \\
5 & 3.131 & 86.3, 0.7189 & 0.007~±~0.002 & 1 & 0.007~±~0.002 & 1 \\
10 & 4.529 & 87.9, 0.7230 & 0.004~±~0.002 & 0 & 0.007~±~0.004 & 1 \\
30 & 9.967 & 93.7, 0.7303 & 0.007~±~0.006 & 11 & 0.012~±~0.008 & 3 \\
100 & 28.608 & 109.4, 0.7358 & 0.010~±~0.024 & 22 & 0.031~±~0.015 & 1 \\
200 & 55.021 & 110.4, 0.7359 & 0.030~±~0.060 & 10 & 0.099~±~0.041 & 0 \\
\end{longtable}

\hypertarget{stochastic-execution-of-the-nominal-schedules}{%
\subsection{6.5. Stochastic execution of the nominal schedules}\label{stochastic-execution-of-the-nominal-schedules}}

The optimisation evaluates routes on their nominal schedule (Assumption 3). To test what that approximation changes under the assumed stochastic model (it does not test the reliability functions themselves, which remain uncalibrated, Section 6.7), every route of the reference Instance 1 fronts (194 hard-dwell, 115 soft-dwell) and of the BBB-20 hard-dwell front (2,045 routes) was executed in simulation, and so was every one of the 12,612,600 Instance 1 routes whose nominal cost lies within 0.5 \% of the optimum for at least one failure cost of Table 8 (506 hard-dwell routes, 452 of them outside the front; 7 soft-dwell routes, all on the front), so that the comparison is not confined to the front. Each travel and processing time is drawn independently from a normal distribution with the model's mean and standard deviation, clipped at zero (negative draws set to zero); the dwell rule is applied in every realisation; and the realised gaps enter the flow reliability, so that the reliability of a realisation is the product of its own pass, flow and travel terms and the expected reliability is their mean. Every route is simulated with its own random stream, seeded by the route itself, so that any two routes are simulated independently and a route that belongs to two candidate sets receives identical draws; Instance 1 routes received 20,000 realisations each and BBB-20 routes 5,000, and the experiment was repeated with both standard deviations doubled. Table 16 compares nominal and simulated values for the key routes. The simulation tests the nominal-schedule approximation on these candidate sets, not over all \(|\Omega|\) routes; statements about selections below refer to the best route among the simulated candidates. The counts of non-dominated routes under simulated objectives are estimates that depend on the sample. Standard errors of expected costs are computed from the per-realisation costs \(c_{t}\widetilde{T} + c_{f}(1 - \widetilde{R})\), so that they include the covariance of time and reliability within a route; a difference between two routes, which are simulated independently, combines the two.

\textbf{Table 16.} Nominal versus simulated values for key Instance 1 routes (20,000 realisations per route; standard errors at most 0.007 s for the mean cycle time and \(10^{-4}\) for the reliability).

\begin{longtable}[]{@{}
  >{\raggedright\arraybackslash}p{(\columnwidth - 8\tabcolsep) * \real{0.2539}}
  >{\raggedright\arraybackslash}p{(\columnwidth - 8\tabcolsep) * \real{0.2364}}
  >{\raggedright\arraybackslash}p{(\columnwidth - 8\tabcolsep) * \real{0.1756}}
  >{\raggedright\arraybackslash}p{(\columnwidth - 8\tabcolsep) * \real{0.1775}}
  >{\raggedright\arraybackslash}p{(\columnwidth - 8\tabcolsep) * \real{0.1566}}@{}}
\toprule\noalign{}
\begin{minipage}[b]{\linewidth}\raggedright
Route
\end{minipage} & \begin{minipage}[b]{\linewidth}\raggedright
Sequence
\end{minipage} & \begin{minipage}[b]{\linewidth}\raggedright
\(T\) / mean realised (s)
\end{minipage} & \begin{minipage}[b]{\linewidth}\raggedright
\(\sqrt{\sigma^{*}}\) / realised SD (s)
\end{minipage} & \begin{minipage}[b]{\linewidth}\raggedright
\(R\) / mean realised
\end{minipage} \\
\midrule\noalign{}
\endhead
\bottomrule\noalign{}
\endlastfoot
IC-by-IC baseline (hard) & AAAAABBBBCCCDDD & 191.81 / 191.82 & 0.849 / 0.845 & 0.7183 / 0.7183 \\
hard, min \(T\) & ABCADBCADBACDBA & 87.55 / 87.59 & 0.934 / 0.608 & 0.7597 / 0.7598 \\
hard, min \(T\) (tie) & ABABDCABDCABDCA & 87.55 / 87.58 & 0.910 / 0.600 & 0.7576 / 0.7576 \\
hard, selected at \(c_{f} = 30\) \$ & ABDACBDCABACDBA & 90.83 / 90.83 & 0.922 / 0.728 & 0.7639 / 0.7639 \\
hard, max \(R\) & ADCBDCAABBAABDC & 114.31 / 114.30 & 0.896 / 0.900 & 0.7706 / 0.7706 \\
soft, selected at \(c_{f} = 2\) \$ & ABCADBACDBACDBA & 75.94 / 75.93 & 0.923 / 0.925 & 0.7068 / 0.7063 \\
soft, max \(R\) & ABDACBADCBADCBA & 78.61 / 78.61 & 0.954 / 0.954 & 0.7200 / 0.7195 \\
\end{longtable}

Three results follow. First, the nominal cycle time and reliability are accurate. Over the front routes the mean realised cycle time exceeds \(T\) by at most 0.16 s on Instance 1 and 0.26 s on BBB-20 (0.45 s and 0.67 s with doubled variability), and over the near-optimal Instance 1 routes by at most 0.21 s (0.51 s), the excess being the hard-dwell maxima of Proposition 1; the mean realised reliability differs from \(R\) by at most 0.001 (0.004 with doubled variability) on every simulated route, because the realised gaps vary by a few percent of \(\tau\). Second, the cost-based selections are largely stable. For each failure cost of Table 8, the route selected from the nominal objectives is the route of least simulated expected cost among the simulated candidates in 16 of the 18 cases (Instance 1 hard and soft, BBB-20 hard), and in the other two (BBB-20 at \(c_{f} = 2\) and 30 \$) it costs 0.0014 ± 0.0002 and less than 0.0001 ± 0.0002 \$/board more than the simulated best; in the second of these the identity of the sample minimiser changes but not, within sampling error, the cost, and the estimated excess costs are the result to read, not the count of unchanged winners. On Instance 1 the candidates include the 452 near-optimal routes outside the hard-dwell front: at every failure cost the nominal selection remains the simulated best, and the best route outside the front is 3.8--59 standard errors above it. With doubled variability 14 of 18 selections remain the simulated best and the largest excess is 0.0035 \$/board (0.15 \% of the cost); in one of the four other cases, Instance 1 hard at \(c_{f} = 2\) \$, the simulated best is a route outside the nominal front, 0.0009 ± 0.0002 \$/board cheaper than the selection. Routes outside the nominal front can also be non-dominated under the simulated objectives: of the 646 simulated Instance 1 hard-dwell routes, 95 are non-dominated in (mean, standard deviation, mean reliability), 31 of them from outside the front. Third, the realised standard deviation never exceeds \(\sqrt{\sigma^{*}}\) beyond sampling error (largest ratio 1.014 on Instance 1 with 20,000 realisations and 1.026 on BBB-20 with 5,000), as Proposition 1 requires, and the bound is tight for soft-dwell routes (ratio 0.99--1.01, Spearman rank correlation 0.99) but conservative for hard-dwell routes. A forced wait restarts the clock at the end of the same IC's dwell and absorbs the variability accumulated since its previous pass, so the realised standard deviation is 0.64--1.03 times the bound, and the nominal and realised variabilities show only a weak negative rank correlation across the hard-dwell fronts (Spearman \(-0.22\) on Instance 1, \(-0.27\) on BBB-20). Under the simulated objectives, and within each front alone, 66 of the 194 Instance 1 hard-dwell front routes remain non-dominated (80 of 115 soft-dwell; 302 of 2,045 on BBB-20 hard-dwell). For hard-dwell operation \(\sigma^{*}\) is therefore a conservative bound to be used as such; a line that needs the variance of the realised hard-dwell cycle time should obtain it by simulating the few routes it shortlists from the front, as done here, which takes seconds per route.

\hypertarget{dispatch-rules-and-waiting-policies}{%
\subsection{6.6. Dispatch rules and waiting policies}\label{dispatch-rules-and-waiting-policies}}

The 54 \% reduction of Section 6.2 is measured against IC-by-IC processing, which no line would keep once interleaving is allowed. Three dispatch rules provide a practical baseline: round-robin (RR), which cycles through the ICs in index order and skips exhausted ones; nearest-next (NN), which moves to the nearest other IC with passes left; and earliest-ready-nearest (ERN), which moves to the IC whose next pass can start earliest, the start being the later of the arrival and the end of that IC's nominal dwell, with ties broken by distance. Each rule was run under both dwell policies and the cheaper result kept; the reference optimum is the cheapest route on the union of both fronts (Table 17).

\textbf{Table 17.} Dispatch rules against the reference optimum: expected cost per board at \(c_{f} = 30\) \$ (\(c_{t} = 0.02\) \$/s), and the excess cost of round-robin over the reference optimum at \(c_{f} = 5\), 30 and 200 \$.

\begin{longtable}[]{@{}
  >{\raggedright\arraybackslash}p{(\columnwidth - 12\tabcolsep) * \real{0.1899}}
  >{\raggedright\arraybackslash}p{(\columnwidth - 12\tabcolsep) * \real{0.1444}}
  >{\raggedright\arraybackslash}p{(\columnwidth - 12\tabcolsep) * \real{0.0725}}
  >{\raggedright\arraybackslash}p{(\columnwidth - 12\tabcolsep) * \real{0.0864}}
  >{\raggedright\arraybackslash}p{(\columnwidth - 12\tabcolsep) * \real{0.0725}}
  >{\raggedright\arraybackslash}p{(\columnwidth - 12\tabcolsep) * \real{0.1438}}
  >{\raggedright\arraybackslash}p{(\columnwidth - 12\tabcolsep) * \real{0.2906}}@{}}
\toprule\noalign{}
\begin{minipage}[b]{\linewidth}\raggedright
Instance
\end{minipage} & \begin{minipage}[b]{\linewidth}\raggedright
IC-by-IC (\$)
\end{minipage} & \begin{minipage}[b]{\linewidth}\raggedright
NN (\$)
\end{minipage} & \begin{minipage}[b]{\linewidth}\raggedright
ERN (\$)
\end{minipage} & \begin{minipage}[b]{\linewidth}\raggedright
RR (\$)
\end{minipage} & \begin{minipage}[b]{\linewidth}\raggedright
Reference (\$)
\end{minipage} & \begin{minipage}[b]{\linewidth}\raggedright
RR excess at 5 / 30 / 200 \$ (\%)
\end{minipage} \\
\midrule\noalign{}
\endhead
\bottomrule\noalign{}
\endlastfoot
Inst-1 & 12.29 & 10.30 & 9.06 & 9.00 & 8.90 & 1.1 / 1.2 / 3.7 \\
Inst-2 & 11.08 & 9.17 & 7.68 & 7.25 & 7.14 & 2.6 / 1.5 / 1.4 \\
Inst-3 & 13.05 & 10.75 & 9.39 & 8.67 & 8.65 & 1.1 / 0.2 / 1.5 \\
BBB-14 & 10.16 & 8.45 & 7.24 & 6.76 & 6.74 & 0.4 / 0.2 / 0.6 \\
BBB-16 & 11.69 & 9.76 & 8.32 & 7.54 & 7.50 & 0.7 / 0.5 / 0.8 \\
BBB-20 & 15.20 & 12.67 & 10.97 & 10.09 & 9.97 & 3.2 / 1.3 / 3.2 \\
\end{longtable}

Round-robin is close to optimal on every instance: its cost is within 3.7 \% of the reference optimum for failure costs from 5 to 200 \$, because cycling through all ICs gives each one nearly the longest gap available between its visits. The one-step-greedy ERN rule is 1.8--11 \% above the optimum at \(c_{f} = 30\) \$ (up to 13 \% at 200 \$), and NN, which returns to nearby ICs too early, 16--30 \% above it at \(c_{f} = 30\) \$. Interleaving therefore delivers most of the gain over IC-by-IC processing, and a simple rule captures most of interleaving. The reference front adds the last 0.2--3.7 \% at the failure costs of Table 17 (0.01--1.8 \$ per board), and, for a line that already uses round-robin, it is the only way to know that the remaining gap is this small; round-robin's position depends on the IC order and on how evenly flow times and pass counts are spread, which is what the reference front checks for each board design.

The reference fronts were computed under the two waiting policies of Assumption 4. A partial-dwell policy, in which the head waits until the gap reaches \(\alpha\tau_{i}\), interpolates between them (\(\alpha = 0\) soft, \(\alpha = 1\) hard). Solving Instance 1 exactly for \(\alpha = 0.25\), 0.5 and 0.75 adds 13 of the 100 distinct \((T,R)\)-non-dominated points of the union of the five fronts, so intermediate waiting does create new trade-offs. The cost-based selection, however, is unchanged at every failure cost of Table 8, and over a grid of \(c_{f}\) from 0.5 to 300 \$ in steps of 0.05 \$ the largest gain from a partial-dwell route is 0.005 \$/board (0.18 \%), at \(c_{f} \approx 4.6\) \$, where the soft- and hard-dwell cost curves cross. Waiting chosen freely for each pass would enlarge the decision space beyond multiset permutations and is left for future work.

\hypertarget{sensitivity-and-boundaries-of-the-study}{%
\subsection{6.7. Sensitivity and boundaries of the study}\label{sensitivity-and-boundaries-of-the-study}}

Table 18 reports Instance 1 under changes of \(r_{0}\), of all flow times and of the travel hazard \(\lambda\). The optimal interleaving is stable under \(r_{0}\): the same max-\(R\) route is best for \(r_{0} = 0.95\) and \(0.99\), and \(r_{0}\) rescales \(R\) without reordering the front. Flow time is the sensitive parameter. Halving \(\tau\) makes the hard- and soft-dwell max-\(R\) routes coincide (75.9 s, \(R = 0.793\)), although the min-\(T\) members still differ (74.3 s hard, 65.3 s soft); doubling it makes the hard-dwell optimum almost twice as long as the soft-dwell one while the soft optimum loses half its reliability, the regime in which the tri-objective choice matters most. The travel hazard is the parameter that ties reliability to the geometry of the route. Over \(\lambda = 0\) to \(0.002\ s^{-1}\) the minimum-\(T\) time (87.6 s) and the soft-dwell max-\(R\) route are unchanged, the latter rescaling from \(R = 0.736\) to 0.704; the hard-dwell max-\(R\) route shortens from 117.5 s to 113.6 s because long waits and detours now cost reliability; the hard-dwell front shrinks from 354 to 142 points as time and reliability become more aligned; and the route selected at \(c_{f} = 30\) \$ moves from 92.3 s (\(\lambda \leq 0.0005\)) to 90.8 s (\(\lambda \geq 0.001\)) while its cost rises from 8.41 to 9.35 \$/board. The ordering conclusions are therefore robust to \(\lambda\), but a line whose moves are markedly riskier should re-select from its own front.

\textbf{Table 18.} Sensitivity of the Instance 1 fronts to \(r_{0}\), to the flow times and to the travel hazard \(\lambda\): (a) max-\(R\) member, (b) min-\(T\) member.

\begin{enumerate}
\def\labelenumi{(\alph{enumi})}
\tightlist
\item
  max-\(R\) member
\end{enumerate}

\begin{longtable}[]{@{}
  >{\raggedright\arraybackslash}p{(\columnwidth - 8\tabcolsep) * \real{0.3184}}
  >{\raggedright\arraybackslash}p{(\columnwidth - 8\tabcolsep) * \real{0.1342}}
  >{\raggedright\arraybackslash}p{(\columnwidth - 8\tabcolsep) * \real{0.3298}}
  >{\raggedright\arraybackslash}p{(\columnwidth - 8\tabcolsep) * \real{0.1266}}
  >{\raggedright\arraybackslash}p{(\columnwidth - 8\tabcolsep) * \real{0.0910}}@{}}
\toprule\noalign{}
\begin{minipage}[b]{\linewidth}\raggedright
Change
\end{minipage} & \begin{minipage}[b]{\linewidth}\raggedright
Variant
\end{minipage} & \begin{minipage}[b]{\linewidth}\raggedright
Route
\end{minipage} & \begin{minipage}[b]{\linewidth}\raggedright
\(T\) (s)
\end{minipage} & \begin{minipage}[b]{\linewidth}\raggedright
\(R\)
\end{minipage} \\
\midrule\noalign{}
\endhead
\bottomrule\noalign{}
\endlastfoot
\(r_{0} = 0.95\) & soft & ABDACBADCBADCBA & 78.6 & 0.524 \\
\(r_{0} = 0.99\) (base) & soft & ABDACBADCBADCBA & 78.6 & 0.720 \\
\(r_{0} = 0.999\) & soft & ABDACBACDBADCBA & 77.5 & 0.782 \\
\(\tau \times 0.5\) & hard & ABCADBACDBACDBA & 75.9 & 0.793 \\
\(\tau \times 0.5\) & soft & ABCADBACDBACDBA & 75.9 & 0.793 \\
\(\tau \times 2.0\) & hard & CDBAAABBDCAABDC & 183.2 & 0.767 \\
\(\tau \times 2.0\) & soft & ABDACBADCBADCBA & 78.6 & 0.350 \\
\(\lambda = 0\) & hard & ADCBDAACBBCDAAB & 117.5 & 0.784 \\
\(\lambda = 0.0005\) & hard & ADCBDAABBDCAABC & 115.0 & 0.777 \\
\(\lambda\) base (0.001) & hard & ADCBDCAABBAABDC & 114.3 & 0.771 \\
\(\lambda = 0.002\) & hard & CADBCDAABBAABDC & 113.6 & 0.758 \\
\(\lambda = 0\) & soft & ABDACBADCBADCBA & 78.6 & 0.736 \\
\(\lambda = 0.002\) & soft & ABDACBADCBADCBA & 78.6 & 0.704 \\
\end{longtable}

\begin{enumerate}
\def\labelenumi{(\alph{enumi})}
\setcounter{enumi}{1}
\tightlist
\item
  min-\(T\) member
\end{enumerate}

\begin{longtable}[]{@{}
  >{\raggedright\arraybackslash}p{(\columnwidth - 8\tabcolsep) * \real{0.4233}}
  >{\raggedright\arraybackslash}p{(\columnwidth - 8\tabcolsep) * \real{0.1135}}
  >{\raggedright\arraybackslash}p{(\columnwidth - 8\tabcolsep) * \real{0.2790}}
  >{\raggedright\arraybackslash}p{(\columnwidth - 8\tabcolsep) * \real{0.1071}}
  >{\raggedright\arraybackslash}p{(\columnwidth - 8\tabcolsep) * \real{0.0770}}@{}}
\toprule\noalign{}
\begin{minipage}[b]{\linewidth}\raggedright
Change
\end{minipage} & \begin{minipage}[b]{\linewidth}\raggedright
Variant
\end{minipage} & \begin{minipage}[b]{\linewidth}\raggedright
Route
\end{minipage} & \begin{minipage}[b]{\linewidth}\raggedright
\(T\) (s)
\end{minipage} & \begin{minipage}[b]{\linewidth}\raggedright
\(R\)
\end{minipage} \\
\midrule\noalign{}
\endhead
\bottomrule\noalign{}
\endlastfoot
\(r_{0} = 0.95\) / 0.99 / 0.999 & soft & CCCDDDBBBBAAAAA & 65.3 & ≈ 0 \\
\(\tau \times 0.5\) & hard & CABADBABDCADBAC & 74.3 & 0.784 \\
\(\tau \times 0.5\) & soft & CCCDDDBBBBAAAAA & 65.3 & ≈ 0 \\
\(\tau \times 2.0\) & hard & ABDCABDCABABDCA & 147.6 & 0.757 \\
\(\tau \times 2.0\) & soft & CCCDDDBBBBAAAAA & 65.3 & ≈ 0 \\
\(\lambda = 0\) & hard & ABCADBCADBADCBA & 87.6 & 0.775 \\
\(\lambda = 0.0005\) & hard & ABCADBCADBACDBA & 87.6 & 0.767 \\
\(\lambda = 0.001\) (base) & hard & ABCADBCADBACDBA & 87.6 & 0.760 \\
\(\lambda = 0.002\) & hard & ABCADBCADBACDBA & 87.6 & 0.745 \\
\end{longtable}

The reference fronts were obtained from completed searches under the stated model, nominal schedules, the two waiting policies of Assumption 4 and parametric reliability functions that were not calibrated on inspection data, and under the numerical conventions of Section 5.2. Within it, the conclusions about route ordering are robust to \(r_{0}\) and \(\lambda\), the nominal selections incur small estimated cost penalties on the simulated candidate sets (Section 6.5), and partial-dwell policies do not change them (Section 6.6); the absolute values of \(R\) are not to be read as field probabilities, and the pass parameters of Section 6.1 are illustrative, which is why \(R/K_{P}\) is reported alongside \(R\). The most valuable further evidence would be one inspected device and recipe family linking the inter-pass gap to an explicitly defined acceptance outcome. NSGA-II and MSSO were run with textbook parameters and no tuning; tuned versions, or methods with explicit diversity maintenance in objective space, may recover more of the interior, and the reference fronts are the ground truth against which such improvements can be measured. The comparison is at equal evaluation budgets, not equal wall time; at \(10^{7}\) evaluations a heuristic run takes seconds where the exact search for BBB-20 takes 38 minutes on 20 cores, and the regret figures of Table 15 are what that difference buys. The exact method is not the fastest conceivable one either: the bound of Equation (11) completes the remaining passes with \(t_{\min}\) alone, whereas the sum of each remaining IC's cheapest outgoing travel would be tighter and still monotone, and a value ordering of \(A_{j}\) by travel time would make the archive tighten earlier. The regime of completed exact search for the present implementation is \(1.5 \times 10^{12}\) routes with a warm start; larger boards may require a tighter bound, a GPU leaf evaluation, or both. The result of Algorithm 1 is determined in exact arithmetic; in the floating-point implementation, the parallel program reproduced each reference front identically across the three campaigns (Table 11) and the matched serial builds returned the same fronts (Table 12). Wall times are reported as mean ± SD over five repeats where repeats were run and as single measurements for the BBB-20 runs and the matched baselines.

\hypertarget{conclusions}{%
\section{7. Conclusions}\label{conclusions}}

The order in which a dispensing head executes the passes of a multi-pass underfill recipe is a routing decision whose reliability is created by the schedule it induces, not carried by the arcs it uses. Formulating it as the TR-UDSP, with a flow reliability that depends on the inter-pass gap, makes the trade between cycle time and void risk explicit as a Pareto front, and the expected-cost or requirement-based selection turns that front into a decision the line engineer can audit. The multiset-permutation BAT enumerates exactly the feasible routes, once each; with a prefix dominance bound, a warm start and counter-window parallelism it completes the exact search for boards of \(10^{9}\) routes in seconds and the \(1.47 \times 10^{12}\) routes of a benchmark derived from a real PCB layout, with assigned underfill recipes, in 38 minutes on 20 cores, with a traversal memory of \(O(|N|)\) states where best-first search exhausts 3 GB on boards ten thousand times smaller. The additive variability measure is the exact cycle-time variance under soft dwell and a proven upper bound on it under hard dwell. Simulated stochastic execution of the front routes, and of every near-optimal route on the smallest board, agrees with the nominal cycle times and reliabilities within 0.3 s and 0.001 at the modelled duration variability on the tested candidate sets; among those candidates the nominally cost-selected routes incur maximum estimated excess costs of 0.0014 \$/board at the modelled variability and 0.0035 \$/board with it doubled; and the hard-dwell variability bound is conservative and should be supplemented by simulation for the few routes a line shortlists. On the instances studied, interleaving cuts cycle time by about half relative to IC-by-IC processing at higher reliability, round-robin dispatch captures most of that gain and comes within 3.7 \% of the reference optimum cost, and NSGA-II and a multi-objective simplified swarm optimisation recover nearly the whole front on boards up to \(4 \times 10^{9}\) routes but only part of it on the largest board (about 30 \% of the hard-dwell and 55--67 \% of the soft-dwell points), while the routes they would select under the cost model have a cost regret of at most 0.53 \%. The exact method therefore does not replace the dispatch rule or the heuristic on the shop floor; it establishes, once per board design, how far their answers can be from the best one. Future work will calibrate the flow-reliability function on void-inspection data, tighten the prefix bound and order the digit choices by travel time, treat waiting as a decision, and extend the model to two-face boards and to adaptive re-sequencing from in-line inspection.

\hypertarget{data-and-code-availability}{%
\section*{Data and Code Availability}\label{data-and-code-availability}}
\addcontentsline{toc}{section}{Data and Code Availability}

The supplementary material contains the MP-BAT sources as run in the campaigns (mpbat\_omp v0.4, mpbat\_serial v0.2) and their corrected versions (v0.5, v0.4), the heuristic sources (moheur v0.1), the CP-SAT control (Python, OR-Tools {[}7{]}) and its results, the six instance files and the regression instance of control (vi), the logs of campaigns 1 and 3 with their repository commits (f1d2322, 30e1226; the log of campaign 2, commit 8bee08b, which produced the heuristic archives, is available from the author), the complete reference fronts as CSV, the 1,440 heuristic archives with their per-run scores, the validation programs (mpbat\_lab v0.2, near v0.1, mc v0.3 with the analysis script that reproduces every statistic of Section 6.5, dispatch v0.1, suffix v0.1) with their outputs, and the earlier simulation versions as records, a script that regenerates Tables 12, 14, 16 and 17 and Figure 3 from those files, a one-command test suite that rebuilds the solvers and runs the regression, encoding and seed checks, and a SHA-256 manifest of every file. The historical campaign launchers are included as records of what was run and are not portable outside the campaign directory layout. All campaign runs were executed on an NVIDIA DGX Spark (host spark-fe05, 20-core Grace CPU, aarch64, 121 GB RAM, Ubuntu 24.04, gcc 13.3, -O3, OMP\_NUM\_THREADS = 20); the CP-SAT control, the single-core development runs and the validation programs used one x86-64 core of an Ubuntu 24.04 container. The BeagleBone Black assembly XY report is available from the BeagleBoard.org repository {[}29{]} under CC BY-SA 3.0.

All six benchmark instances are public: Instances 1--3 are specified completely in Section 6.1, and the three BBB instances are derived from the BeagleBone Black Rev C assembly XY report {[}29{]}, whose SHA-256 digest (68483dd5\ldots c34e9c28) was verified against an independent download, with the pass recipes given in the instance files.

\hypertarget{funding}{%
\section*{Funding}\label{funding}}
\addcontentsline{toc}{section}{Funding}

This research was supported in part by the National Science and Technology Council (NSTC) of Taiwan under grants NSTC 114-2923-E-007-003-MY2 and NSTC 114-2918-I-007-001, and by Wistron NeWeb Corporation (WNC), ``Development of a prediction and optimisation system for advanced-packaging dispensing processes'' (project title translated from Chinese), under Industry--Academia Collaboration Project 114A0271JB. Research supported by the NVIDIA Academic Grant Program (project ``AMR Deployment via RL, Sensor Fusion \& SSO Optimization'') through two RTX PRO 6000 Blackwell Max-Q workstation GPUs for the computing infrastructure of the Integrated \& Collaborative AI Robotics Laboratory (ICARL); the results reported here were computed on CPU cores only (Section 6.1).

\hypertarget{declaration-of-competing-interest}{%
\section*{Declaration of Competing Interest}\label{declaration-of-competing-interest}}
\addcontentsline{toc}{section}{Declaration of Competing Interest}

The author declares no known competing financial interests or personal relationships that could have appeared to influence the work reported in this paper.

\hypertarget{declaration-of-generative-ai-and-ai-assisted-technologies}{%
\section*{Declaration of Generative AI and AI-Assisted Technologies}\label{declaration-of-generative-ai-and-ai-assisted-technologies}}
\addcontentsline{toc}{section}{Declaration of Generative AI and AI-Assisted Technologies}

During the preparation of this work the author used Claude (Anthropic) for drafting assistance, code generation for the enumeration experiments and post-experiment auditing of the manuscript text. The author reviewed and edited the content and takes full responsibility for the content of the publication.

\hypertarget{references}{%
\section*{References}\label{references}}
\addcontentsline{toc}{section}{References}

{[}1{]} Zhang Z, Wong CP. Recent advances in flip-chip underfill: materials, process, and reliability. IEEE Trans Adv Packag 2004;27(3):515--24. \url{https://doi.org/10.1109/TADVP.2004.831870}.

{[}2{]} Wan JW, Zhang WJ, Bergstrom DJ. Recent advances in modeling the underfill process in flip-chip packaging. Microelectron J 2007;38(1):67--75. \url{https://doi.org/10.1016/j.mejo.2006.09.017}.

{[}3{]} Young WB, Yang WL. Underfill viscous flow between parallel plates and solder bumps. IEEE Trans Compon Packag Technol 2002;25(4):695--700. \url{https://doi.org/10.1109/TCAPT.2002.806176}.

{[}4{]} Lewis A. Underfill techniques: automated dispensing and jetting. Solid State Technology (Semiconductor Digest), June 2003. \url{https://sst.semiconductor-digest.com/2003/06/underfill-techniques-automated-dispensing-and-jetting/} (accessed 25 September 2026).

{[}5{]} Hung HH, Cheng YC, Hwang SJ, Chen DL, Chang HJ, Huang BY, et al.~Analysis of flip-chip ball grid array underfill flow process. Int J Adv Manuf Technol 2024;134(9):4851--70. \url{https://doi.org/10.1007/s00170-024-14304-1}.

{[}6{]} Yeh WC. Time-reliability optimization for the stochastic traveling salesman problem. Reliab Eng Syst Saf 2024;248:110179. \url{https://doi.org/10.1016/j.ress.2024.110179}.

{[}7{]} Perron L, Didier F. The CP-SAT-LP solver (invited talk). In: 29th International Conference on Principles and Practice of Constraint Programming (CP 2023). LIPIcs, vol.~280. Dagstuhl: Schloss Dagstuhl--Leibniz-Zentrum für Informatik; 2023, p.~3:1--3:2. \url{https://doi.org/10.4230/LIPIcs.CP.2023.3}.

{[}8{]} Abas A, Ng FC, Gan ZL, Ishak MHH, Abdullah MZ, Chong GY. Effect of scale size, orientation type and dispensing method on void formation in the CUF encapsulation of BGA. Sādhanā 2018;43:59. \url{https://doi.org/10.1007/s12046-018-0849-3}.

{[}9{]} Lin YK, Chang PC. Evaluate the system reliability for a manufacturing network with reworking actions. Reliab Eng Syst Saf 2012;106:127--37. \url{https://doi.org/10.1016/j.ress.2012.05.011}.

{[}10{]} Lin YK, Huang CF, Chang PC. System reliability evaluation of a touch panel manufacturing system with defect rate and reworking. Reliab Eng Syst Saf 2013;118:51--60. \url{https://doi.org/10.1016/j.ress.2013.04.007}.

{[}11{]} Chang PC, Lin YK, Chiang YM. System reliability estimation and sensitivity analysis for multi-state manufacturing network with joint buffers---A simulation approach. Reliab Eng Syst Saf 2019;188:103--9. \url{https://doi.org/10.1016/j.ress.2019.03.024}.

{[}12{]} Lin YK, Huang CF. Reliability evaluation according to a routing scheme for multi-state computer networks under assured accuracy rate. Ann Oper Res 2016;244(1):221--40. \url{https://doi.org/10.1007/s10479-014-1673-8}.

{[}13{]} Niu YF, Song YF, Xu XZ, Zhao X. Efficient reliability computation of a multi-state flow network with cost constraint. Reliab Eng Syst Saf 2022;222:108393. \url{https://doi.org/10.1016/j.ress.2022.108393}.

{[}14{]} Escudero LF. An inexact algorithm for the sequential ordering problem. Eur J Oper Res 1988;37(2):236--49. \url{https://doi.org/10.1016/0377-2217(88)90333-5}.

{[}15{]} Brucker P, Knust S. Complex scheduling. 2nd ed.~Berlin, Heidelberg: Springer; 2012. \url{https://doi.org/10.1007/978-3-642-23929-8}.

{[}16{]} Pop PC, Cosma O, Sabo C, Sitar CP. A comprehensive survey on the generalized traveling salesman problem. Eur J Oper Res 2024;314(3):819--35. \url{https://doi.org/10.1016/j.ejor.2023.07.022}.

{[}17{]} Deb K, Pratap A, Agarwal S, Meyarivan T. A fast and elitist multiobjective genetic algorithm: NSGA-II. IEEE Trans Evol Comput 2002;6(2):182--97. \url{https://doi.org/10.1109/4235.996017}.

{[}18{]} Yeh WC, Zhu W, Yin Y, Huang CL. Cloud computing considering both energy and time solved by two-objective simplified swarm optimization. Appl Sci 2023;13(4):2077. \url{https://doi.org/10.3390/app13042077}.

{[}19{]} Yeh WC, Su YZ, Gao XZ, Hu CF, Wang J, Huang CL. Simplified swarm optimization for bi-objective active reliability redundancy allocation problems. Appl Soft Comput 2021;106:107321. \url{https://doi.org/10.1016/j.asoc.2021.107321}.

{[}20{]} Yeh WC, Zhu W, Tan SY, Wang GG, Yeh YH. Novel general active reliability redundancy allocation problems and algorithm. Reliab Eng Syst Saf 2022;218:108167. \url{https://doi.org/10.1016/j.ress.2021.108167}.

{[}21{]} Yeh WC, Lin YP, Liang YC, Lai CM, Huang CL. Simplified swarm optimization for hyperparameters of convolutional neural networks. Comput Ind Eng 2023;177:109076. \url{https://doi.org/10.1016/j.cie.2023.109076}.

{[}22{]} Yeh WC. Novel binary-addition tree algorithm (BAT) for binary-state network reliability problem. Reliab Eng Syst Saf 2021;208:107448. \url{https://doi.org/10.1016/j.ress.2021.107448}.

{[}23{]} Yeh WC. A quick BAT for evaluating the reliability of binary-state networks. Reliab Eng Syst Saf 2021;216:107917. \url{https://doi.org/10.1016/j.ress.2021.107917}.

{[}24{]} Ko CW, Ruskey F. Generating permutations of a bag by interchanges. Inf Process Lett 1992;41(5):263--9. \url{https://doi.org/10.1016/0020-0190(92)90170-Z}.

{[}25{]} Knuth DE. The art of computer programming, vol.~4A: combinatorial algorithms, part 1. Upper Saddle River (NJ): Addison-Wesley; 2011. ISBN 978-0-201-03804-0.

{[}26{]} Efron B, Stein C. The jackknife estimate of variance. Ann Stat 1981;9(3):586--96. \url{https://doi.org/10.1214/aos/1176345462}.

{[}27{]} Steele JM. An Efron--Stein inequality for nonsymmetric statistics. Ann Stat 1986;14(2):753--8. \url{https://doi.org/10.1214/aos/1176349952}.

{[}28{]} Held M, Karp RM. A dynamic programming approach to sequencing problems. J Soc Ind Appl Math 1962;10(1):196--210. \url{https://doi.org/10.1137/0110015}.

{[}29{]} BeagleBoard.org Foundation. BeagleBone Black Rev C design files and assembly XY report. \url{https://github.com/beagleboard/beaglebone-black} (accessed 15 September 2026).

{[}30{]} Zitzler E, Thiele L. Multiobjective evolutionary algorithms: a comparative case study and the strength Pareto approach. IEEE Trans Evol Comput 1999;3(4):257--71. \url{https://doi.org/10.1109/4235.797969}.

{[}31{]} Coello Coello CA, Reyes Sierra M. A study of the parallelization of a coevolutionary multi-objective evolutionary algorithm. In: MICAI 2004: Advances in Artificial Intelligence. Lecture Notes in Computer Science, vol.~2972. Berlin, Heidelberg: Springer; 2004, p.~688--97. \url{https://doi.org/10.1007/978-3-540-24694-7_71}.

{[}32{]} Zio E, Bazzo R. A clustering procedure for reducing the number of representative solutions in the Pareto Front of multiobjective optimization problems. Eur J Oper Res 2011;210(3):624--34. \url{https://doi.org/10.1016/j.ejor.2010.10.021}.

\end{document}